\documentclass[10pt]{article}

\usepackage{amsmath}
\usepackage{latexsym}
\usepackage{amsfonts}
\usepackage{amssymb}
\usepackage{amscd}
\usepackage{oldgerm}		

\def\define{\def}

\define\M{{\cal{M}}}
\define\Y{{\cal{Y}}}
\define\X{{\cal{X}}}
\define\V{{\cal{V}}}
\define\U{{\cal{U}}}
\define\F{{\cal{F}}}
\define\L{{\cal{L}}}
\define\calS{{\cal{S}}}
\define\I{{\cal{I}}}
\define\E{{\cal{E}}}
\define\W{{\cal{W}}}
\define\D{{\cal{D}}}
\define\O{{\cal{O}}}

\define\C{\Bbb{C}}
\define\R{\Bbb{R}}
\define\Q{\Bbb{Q}}
\define\Z{\Bbb{Z}}

\define\Gc{G_{\C}}
\define\gg{{\textswab{g}}}	
\define\half{\frac{1}{2}}
\define\Hk{H^1(M\setm\Bbb{S},\C)}
\define\hph{\hphantom}
\define\setm{\setminus}

\define\Res{\operatornamewithlimits{Res}}
\define\Im{\operatornamewithlimits{Im}}
\define\Re{\operatorname{Re}}
\define\pd{\partial}
\define\rel{{}^r}

\define\a{\alpha}
\define\b{\beta}
\define\g{\gamma}
\define\l{\lambda}

\define\lam{\Lambda^{-1,-1}}
\define\nn{\eta_{-}}
\define\np{\eta_{+}}
\define\nz{\eta_0}
\define\pip{\pi_{+}}
\define\pin{\pi_{-}}
\define\piz{\pi_0}
\define\pil{\pi_{\Lambda}}

\newtheorem{theorem}{Theorem}[section]
\newtheorem{lemma}{Lemma}[section]
\newtheorem{definition}{Definition}[section]

\newtheorem{corollary}{Corollary}[section]

\newtheorem{example}{Example}[section]

\newtheorem{myremark}{Remark.} 

\numberwithin{equation}{section}

\newbox\qedbox
\setbox\qedbox=\hbox{}

\newenvironment{proof}{\smallskip\noindent{\bf Proof.}\hskip \labelsep}
                        {\hfill\penalty10000\copy\qedbox\par\medskip}

\begin{document}

\begin{center}{\bf Variations of Mixed Hodge Structure, Higgs Fields, 
       and \newline Quantum Cohomology .}
\end{center}
\vskip .5cm

\centerline{Gregory J Pearlstein}
\vskip 1cm


\centerline{\bf Abstract} Following C. Simpson, we show that every variation 
of graded-polarized mixed Hodge structure defined over $\Q$ carries a natural 
Higgs bundle structure $\bar\pd+\theta$ which is invariant under the $\C^*$ 
action studied in \cite{Simpson}.  We then specialize our construction to 
the context of \cite{Deligne}, and show that the resulting Higgs field 
$\theta$ determines (and is determined by) the Gromov-Witten potential of the 
underlying family of Calabi--Yau threefolds.
\vskip 1cm

\section{Introduction}
\label{intro}

\define\H{{\cal{H}}}

\par Let $X$ be a compact K\"ahler manifold.  Then, by virtue of 
\cite{Simpson}, it is known that a monodromy representation
\begin{equation}
	\rho:\pi_1(X,x_0)\to GL(V)			
\end{equation}
arises from a variation of pure, polarized Hodge structure $\V\to X$ if and 
only if 
\begin{description}
\item[(1)] $\rho$ is semisimple.
\item[(2)] The associated Higgs bundle structure $\bar\pd+\theta$ is invariant
	   under the $\C^*$ action $(1.9)$.
\item[(3)] The representation $\rho$ is defined over $\R$ relative to some
	   choice of real structure $V_{\R}$ on $V$.
\end{description}

\par To study the non-semisimple representations of $\pi_1(X,x_0)$, we begin 
with a variation of graded-polarized mixed Hodge structure
$$
	\V\to X
$$ 
and ask whether the underlying $C^{\infty}$ vector bundle of $\V$ carries a
natural Higgs bundle structure $\bar\pd+\theta$ invariant under the
$\C^*$ action $(1.9)$ studied in \cite{Simpson}.  Our main result is as 
follows:

\setcounter{section}{5}
\setcounter{theorem}{0}					      
\begin{theorem} Let $V\to S$ be a variation of graded-polarized mixed Hodge 
structure, and $\{\U^p\}$ denote the collection of $C^{\infty}$ subbundles 
of $\V$ defined by the rule:
$$
	\U^p_s = \bigoplus_q\, I^{p,q}_{(\F_s,\W_s)}
$$
Then, relative to the Gauss--Manin connection $\nabla$, the direct sum 
decomposition
$$
	\V = \bigoplus_p\,\U^p		
$$
defines a (unpolarized) complex variation of Hodge structure.
\end{theorem}
\setcounter{section}{1}

\par Therefore, by virtue of \cite{Simpson}, the underlying $C^{\infty}$ 
vector bundle $E$ of $\V$ does indeed carry a natural Higgs bundle structure
$\bar\pd+\theta$.  Moreover, because the Higgs bundle structure 
$\bar\pd+\theta$ produced by Theorem $(5.1)$ arises from a complex variation 
of Hodge structure, it is automatically a fixed point of $(1.9)$.
[See Lemma $(1.12)$ for details.]

\begin{myremark} In the case of variations of graded-polarized mixed Hodge 
structure arising from the cohomology of a family of singular or 
quasi-projective varieties, the corresponding Higgs field $\theta$ produced 
by Theorem $(5.1)$ turns out to be a natural analog of the Kodaira--Spencer 
map associated to a smooth family of non-singular projective varieties. 
\end{myremark}

\par The main tools used in proof of Theorem $(5.1)$ are Deligne's observation
that every mixed Hodge structure $(F,W)$ determines a functorial bigrading
$$
	V = \bigoplus_{p,q}\,I^{p,q}_{(F,W)}		
$$
of the underlying vector space $V$ which mimics the classical Hodge 
decomposition (although in general $\bar I^{p,q} \neq I^{q,p}$), and
the classifying spaces $\M$ of graded-polarized mixed Hodge structure
described in \S 3.

\begin{myremark} An alternative proof of Theorem $(5.1)$ due to P. Deligne
\cite{Deligne2} is outlined in the appendix.
\end{myremark}

\par The general outline of this paper is as follows: After presenting some 
basic definitions at the end of this section and the preliminary remarks of 
\S 2, we begin \S 3 with a review of the classifying spaces $\D$ of pure, 
polarized Hodge structure constructed in \cite{Griffiths2}.  Following 
\cite{Kaplan}, we then define classifying space of graded-polarized mixed 
Hodge structure which are universal with respect to variations of 
graded-polarized mixed Hodge structure. 


\par In analogy with the pure case, the fundamental $C^{\infty}$ vector 
bundles supported by such classifying spaces are the turn out to be the 
Deligne--Hodge bundles 
$$
	\I^{p,q}_F = I^{p,q}_{(F,W)}
$$
These bundles are studied in \S 4, wherein we extend the methods of 
\cite{Deligne0} to obtain an explicit formula governing the first 
order behavior of the decomposition:
$$
	V = \bigoplus_{p,q}\,\I^{p,q}_{(F,W)}
$$
along $\M$.  Careful application of this formula then gives the proof of
Theorem $(5.1)$ presented at the beginning of \S 5.
\vskip 5pt

\par To study the asymptotic behavior of such Higgs bundles, we then prove a
\lq\lq group theoretic\rq\rq{} version of Schmid's Nilpotent Orbit Theorem 
[Formula $(6.8)$] and derive an equivalence of categories theorem 
for unipotent variations of mixed Hodge structure [Theorem $(6.19)$].
\vskip 5pt

\par Armed with these results, we then specialize our constructions to the
context of mirror symmetry and quantum cohomology in \S 7 and \S 8.  In
doing so, we obtain interpretations of both Deligne's work on the local 
behavior of Hodge structures at infinity \cite{Deligne} and some recent 
work of David Cox and Sheldon Katz \cite{Cox--Katz} in terms of Higgs 
fields associated to variations of graded-polarized mixed Hodge structures.
\vskip 5pt


\par To close this section, we shall now recall several basic definitions
and constructions which are used throughout this paper.  Additional background
material may be found in \S 2.

\setcounter{definition}{1}					
\begin{definition} Let $S$ be a complex manifold.  Then, following 
\cite{S--Z}, we define a variation $\V\to S$ of graded-polarized mixed Hodge 
structure to consist of a $\Q$-local system $\V_{\Q}$ defined over $S$
equipped with:
\begin{description}
\item[(1)] A rational, increasing weight filtration 
	   $$
		0\subseteq\cdots\W_k\subseteq\W_{k+1}\subseteq\cdots
		 \subseteq\V_{\C}
	   $$
	   of $\V_{\C} = \V_{\Q}\otimes\C$.
\item[(2)] A decreasing Hodge filtration 
	   $$
		0\subseteq\cdots\F^p\subseteq\F^{p-1}\subseteq\cdots
		 \subseteq\V_{\C}\otimes\O_S
	   $$
\item[(3)] A collection of rational, non-degenerate bilinear forms
	   $$
		\calS_k:Gr^{\W}_k(\V_{\Q})\otimes Gr^{\W}_k(\V_{\Q})\to\Q
	   $$
	   of alternating parity $(-1)^k$.
\end{description} 
satisfying the following mutual compatibility conditions:
\begin{description}
\item[(a)] Relative to the Gauss--Manin connection of $\V$:
	   $$
		\nabla\,\F^p\subseteq\Omega^1_S\otimes\F^{p-1}
	   $$
	   for each index $p$.
\item[(b)] The triple $(Gr^{\W}_k(\V_{\Q}),\F\,Gr^{\W}_k(\V_{\Q}),\calS_k)$
	   defines a variation of pure, polarized Hodge structure for each
	   index $k$.
\end{description}
\end{definition}

\setcounter{definition}{2}					
\begin{definition} A Higgs bundle $(E,\bar\pd+\theta)$ consists of a 
holomorphic vector bundle $(E,\bar\pd)$ endowed with a endomorphism 
valued 1-form 
$$
	\theta:\E^0(E)\to \E^{1,0}(E)
$$ 
which is both holomorphic and symmetric (i.e. $\bar\pd\theta = 0$ and 
$\theta\wedge\theta = 0$).  
\end{definition}

\setcounter{example}{3}						
\begin{example} Let $\V$ denote a variation of pure, polarized Hodge 
structure arising via the cohomology of a smooth family of non-singular 
projective varieties $f:Y\to X$.  Then, by virtue of the $C^{\infty}$ 
decomposition
\setcounter{equation}{4}
\begin{equation}
	  \V = \bigoplus_{p+q=k}\,\H^{p,q}			
\end{equation}
underlying smooth vector bundle
$$
	E = \V_{\C}\otimes\E^0_X
$$	
of $\V$ inherits an integrable complex structure $\bar\pd$ via the 
isomorphism
$$
	\H^{p,q} \cong \F^p/\F^{p+1}
$$
and the holomorphic structure of $\F^p$.  Likewise, the Kodaira--Spencer map
$$
	\kappa_p:T_p(X)\to H^1(Y_p,\Theta(Y_p))
$$
defines a symmetric, endomorphism valued 1-form $\theta$ on $E$ via
the rule
$$
	\theta(\xi)(\sigma) = \kappa(\xi)\cup\sigma
$$

\par To prove that $(E,\bar\pd+\theta)$ is indeed a Higgs bundle, observe
that by virtue of \cite{Griffiths2}, we may write the Gauss--Manin 
connection $\nabla$ as
\begin{equation}
	\nabla = \tau+ \bar\pd + \pd +\theta 			
\end{equation}
relative to a pair of differential operators
\begin{equation}
\bar\pd:\E^0(E)\to\E^{0,1}(E),\qquad \pd:\E^0(E)\to\E^{1,0}(E) 
\end{equation}
preserving the Hodge decomposition $(1.5)$ and a pair of tensor fields
\begin{equation}
	\begin{array}{l}
	   \tau:\H^{p,q}\to\E^{0,1}\otimes\H^{p+1,q-1}	  \\
	   \theta:\H^{p,q}\to\E^{1,0}\otimes\H^{p-1,q+1} 
	\end{array}						
\end{equation}
shifting the indices of $(1.5)$ by $\pm 1$.  Expanding out the 
integrability condition $\nabla^2=0$, and taking account of equations
$(1.6)$--$(1.8)$, it then follows that
$$
	\bar\pd^2 = 0,\qquad \bar\pd\theta = 0,\qquad \theta\wedge\theta = 0
$$
and hence $(E,\bar\pd+\theta)$ is a Higgs bundle by virtue of the 
Newlander-Nirenburg theorem.  Moreover, given any element $\l\in\C^*$, the
map
$$
	f:\V\to\V,\qquad \left. f\right|_{\H^{p,q}} = \l^p
$$
defines a bundle isomorphism 
$$
	(E,\bar\pd+\theta) \cong (E,\bar\pd+\l\theta)
$$
Consequently, the isomorphism class of such a Higgs bundle $(E,\bar\pd+\theta)$
is a fixed point of the $\C^*$ action 
\begin{equation}
	\l:(E,\bar\pd+\theta)\to (E,\bar\pd+\l\theta)		
\end{equation}
\end{example}

\par In order to construct deformations of such Higgs bundles, Simpson
introduces the following definition:

\setcounter{definition}{9}					
\begin{definition} A complex variation of Hodge structure consists of the
following data:
\begin{description}
\item[(1)] A flat, $\C$-vector bundle $(E,\nabla)$.
\item[(2)] A $C^{\infty}$ decomposition 
	   \setcounter{equation}{10}
	   \begin{equation}
        	E = \bigoplus_p\,{\cal U}^p			
	   \end{equation}
	   satisfying Griffiths' horizontality, i.e. 
	   $$
	      \nabla:\E^0(\U^p) \to 
		     \E^{0,1}(\U^{p+1})\oplus\E^{0,1}(\U^p)\oplus
		     \E^{1,0}(\U^p)\oplus\E^{1,0}(\U^{p-1})
           $$				
\end{description}
\end{definition}

\begin{myremark} A flat hermitian form $Q$ is said to polarize the complex
variation $(E,\nabla,\oplus_p\,\U^p)$ provided bilinear form
$$
	\langle u,v\rangle = Q(Cu,v),\qquad \left.C\right|_{\U^p} = (-1)^p
$$
is positive definite and makes the direct sum decomposition $(1.11)$ In 
general, the complex variations of Hodge structure considered in this paper 
will be unpolarized.
\end{myremark}

\setcounter{lemma}{11}					
\begin{lemma} Every complex variation of Hodge structure $(E,\nabla,\U^{\ast})$
carries a natural Higgs bundle structure $\bar\pd+\theta$ invariant under the
$\C^*$ action $(1.9)$.
\end{lemma}
\begin{proof} One simply goes through the proof presented in Example $(1.4)$,
replacing $\H^{p,q}$ with $\U^p$.
\end{proof}

$$
	\text{\bf Acknowledgements}
$$

\par The author thanks his thesis advisor Aroldo Kaplan for his patience and 
his guidance.  The author also thanks P. Deligne for his many helpful comments,
and suggestions as well as David Cox and Sheldon Katz for their help in 
connection with the material of \S 7, and Eduardo Cattani for clarifying the
role of Higgs fields in \S 6.  At a much more general level, I must also thank
Ivan Mirkovic, Andrey Todorov and Steven Zucker along with my fellow students
Eric Gaze, Christine Schwarz and Glenn Fawcett for their help, support and
encouragement.

\section{Preliminary Remarks}
\label{sec:2}

\par The purpose of this section is to acquaint the reader with P. Deligne's
theory of mixed Hodge structures, and provide a catalog of basic definitions
for later use.  We assume only that the reader is already familiar with the
basic tenets of Hodge theory, as outlined in \cite{GH}.

\par  Conceptually, a the notion of a mixed Hodge structure may be viewed as
a kind of \lq\lq iterated extension\rq\rq{} of pure Hodge structures.  Thus,
as a prelude to the formal definition presented below, let us first consider 
the problem of defining what it should mean for an exact sequence 
\setcounter{equation}{0}
\begin{equation}
	0 \rightarrow A \stackrel{\a}{\rightarrow} B 
	  \stackrel{\b}{\rightarrow} C \rightarrow 0	
\end{equation}
to define an extension of pure Hodge structures.  Since both $A$ and $C$ are 
filtered vector spaces, one might be tempted to require only that the vector 
space $B=B_{\Q}\otimes\C$ carry a decreasing \lq\lq Hodge filtration\rq\rq{} 
$F^{\bullet}(B)$ which is strictly compatible with the given maps $\a$ and 
$\b$, i.e.
\begin{equation}
	\a(F^p(A)) = F^p(B)\cap\a(A),\quad \b(F^p(B)) = F^p(C)\cap\b(B) 
\end{equation}

\par The trouble with this preliminary definition is that it does not encode 
the weights of the pure Hodge structures $A$ and $C$.  To rectify this 
defect, observe that a pure Hodge structure of weight $k$ may be completely 
recovered from the knowledge of its Hodge filtration $F$ and its 
weight filtration
$$
	W_j(V_{\Q}) = \left\{ \begin{array}{ll}
				V_{\Q} & \qquad j\geq k \\
				0      & \qquad j<k
			      \end{array}
		      \right.
$$
via the rule $H^{p,q} = F^p\cap\bar F^q\cap W_{p+q}$.  In light
of this observation, it is therefore natural to require that $B$ carry 
an increasing \lq\lq weight filtration\rq\rq{}
$$
	0\subseteq\cdots\subseteq W_{k-1}(B_{\Q})\subseteq W_k(B_{\Q})
	 \subseteq\cdots\subseteq B_{\Q}
$$
which is strictly compatible with the given maps $\a$ and $\b$, i.e.
\begin{equation}
	\begin{array}{l}
	\a(W_k(A_{\Q})) = W_k(B_{\Q})\cap\a(A_{\Q}) \\     
	\b(W_k(B_{\Q})) = W_k(C_{\Q})\cap\b(B_{\Q}) 
	\end{array} 				      
\end{equation}

\par To see that these two requirements define a reasonable notion
of what it should mean for $B$ to represent an extension of $C$ by $A$, 
observe that one at least has the following result:

\setcounter{lemma}{3}					
\begin{lemma} Let $F^{\bullet}(B)$ and $W_{\bullet}(B_{\Q})$ be a pair of 
filtrations which satisfy condition $(2.2)$ and $(2.3)$.  Then, for each 
index $k$, the \lq\lq Hodge filtration\rq\rq{} $F^{\bullet}(B)$ 
induces a pure Hodge structure of weight $k$ on the quotient space
$Gr^{W(B)}_k := W_k(B)/W_{k-1}(B)$ via the rule:
$$
	F^p Gr^{W(B)}_k = \frac{F^p(B)\cap W_k(B) + W_{k-1}(B)}{W_{k-1}(B)}
$$
\end{lemma}

\par The pair $(F(B),W(B))$ constructed above is an example of a mixed Hodge
structure.  The geometric importance such structures rests upon P. Deligne's 
construction of functorial mixed Hodge structures on the cohomology of an 
arbitrary complex algebraic variety $X$.

\setcounter{definition}{4}				
\begin{definition} Let $V=V_{\Q}\otimes\C$ be a finite dimensional complex
vector space defined $\Q$.  Then, a decreasing filtration $F$ of $V$ is said 
to pair with an increasing filtration $W=W(V_{\Q})\otimes\C$ to define a mixed 
Hodge structure $(F,W)$ provided they satisfy the following condition:  For 
each index $k$, the induced filtration
$$
	F^p Gr^W_k = \frac{F^p\cap W_k + W_{k-1}}{W_{k-1}}
$$
defines a pure Hodge structure of weight $k$ on quotient space 
$$
	Gr^W_k = W_k/W_{k-1}
$$
\end{definition}

\setcounter{example}{5}					
\begin{example} Let $\Bbb S$ be a finite set of distinct points
in a compact Riemann surface $M$ and $\Omega^1_M(\Bbb S)$ denote the 
space of meromorphic 1-forms on $M$ which have at worst simple poles
along $\Bbb S$.  Then, the mixed Hodge structure $(F,W)$ attached to $\Hk$ 
by Deligne's construction is given by the following pair of filtrations:
\begin{equation*}
  \begin{array}{lll}
   \qquad W_0 = 0     &\qquad W_1 = H^1(M,\C) &\qquad W_2 = \Hk \\
   \qquad F^2 = 0     &\qquad F^1 = \Omega^1_M(\Bbb S)   &\qquad F^0 = \Hk 
\end{array}
\end{equation*}
\end{example}

\par To obtain an analog of the Hodge decomposition in the category of 
mixed Hodge structures, Deligne proceeds as follows:

\setcounter{definition}{6}				
\begin{definition} A bigrading of a mixed Hodge structure $(F,W)$ is a
direct sum decomposition $V = \bigoplus_{p,q}\,J^{p,q}$ of the the 
underlying complex vector space $V$ which has the following  two 
properties:
\begin{description}
\item{$\bullet$} $F^p = \oplus_{r\geq p,s}\, J^{r,s}$.
\item{$\bullet$} $W_k = \oplus_{r+s\leq k}\, J^{r,s}$.
\end{description}
\end{definition}

\setcounter{lemma}{7}					
\begin{lemma} Let $(F,W)$ be a mixed Hodge structure.  Then, there exists a 
unique bigrading $\{I^{p,q}\}$ of $(F,W)$ with the following additional
property:
$$
		I^{p,q} = \bar I^{q,p} \mod \bigoplus_{r<p,s<q}\, I^{r,s}
$$
\end{lemma}

\setcounter{example}{8}					
\begin{example} In the case of the finitely punctured Riemann surface 
$M\setm\Bbb S$ considered previously, the {\it Deligne--Hodge decomposition}
of $(F,W)$ is given by the following subspaces of $\Hk$:
$$
    I^{1,1} = F^1\cap\overline{F^1},\qquad I^{1,0} = H^{1,0}(M),
    \qquad I^{0,1} = H^{0,1}(M)
$$
Moreover, in this particular case, the subspace $I^{1,1}\subseteq F^1$
admits the following description:  Let ${\cal H}^0_M(\Bbb S)$ denote the space
of real-valued harmonic functions on $M$ which have at worst logarithmic 
singularities along $\Bbb S$.  Then,
$$
	I^{1,1}\cap H^1(M\setminus\Bbb S,\R)
	= \left\{\, \sqrt{-1}\,\frac{\pd f}{\pd z}dz \mid 
			f\in {\cal H}^0_M(\Bbb S) \,\right\}
$$
\end{example}

\begin{myremark}{\bf [Kaplan]} By virtue of Lemma $(2.8)$, the subspaces
$$
	\Lambda^{p,q}(V) := \bigoplus_{a\leq p, b\leq q}\, I^{p,q}
$$
satisfy the symmetry condition $\overline{\Lambda^{p,q}(V)}=\Lambda^{q,p}(V)$.
In particular, 
\setcounter{equation}{9}					
\begin{equation}
	\bar I^{p,q} = I^{q,p} \mod \Lambda^{q-1,p-1}(V) 
\end{equation}
\end{myremark}

\par Regarding the functorial properties of mixed Hodge structures, one has
the following basic result:

\setcounter{theorem}{10}			   
\begin{theorem} The category of mixed Hodge structures defined over a fixed
subfield $R\subseteq\R$ is abelian.  Moreover, it is closed under the 
operations of taking direct sums, tensor products, and duals.
\end{theorem}

\setcounter{corollary}{11}			    
\begin{corollary} The choice of a mixed Hodge structure $(F,W)$ on 
$V=V_{\Q}\otimes\C$ induces a mixed Hodge structure on $gl(V)$ via the 
bigrading:
$$
	gl(V)^{r,s} = \{\, \a\in gl(V) \mid \a:I^{p,q}\to I^{p+r,q+s}\quad 
					    \forall\,p,q  \,\}
$$
\end{corollary}

\begin{myremark} A morphism of mixed Hodge structure $f:V\to V'$ is a
$R$-linear map which is strictly compatible with the filtrations
$F$ and $W$.
\end{myremark}


\setcounter{definition}{12}				
\begin{definition} A graded-polarization of a mixed Hodge structure 
$(F,W)$ consists of a choice of polarization ${\cal S}_k$ for each 
non-trivial layer $Gr^W_k$ of $Gr^W$.
\end{definition}

\setcounter{example}{13}				
\begin{example} Given a finitely puncture Riemann surface $M\setm\Bbb S$, 
the mixed Hodge structures attached to $\Hk$ is graded-polarized by
the nondegenerate bilinear forms 
$$ {\cal} S_2(\alpha,\beta) =
     4\pi^2\sum_{p\in\Bbb S}\,\text{Res}_p(\alpha)\text{Res}_p(\beta),\qquad
   {\cal} S_1(\alpha,\beta) = \int_M\, \alpha\wedge\beta
$$
defined on $Gr^W_2$ and $Gr^W_1$ respectively.
\end{example}

\begin{myremark} Let $\cal S$ be a graded-polarization of the mixed Hodge 
structure $(F,W)$ and 
$$
	\Gc = \{\, g\in GL(V)^W \mid Gr(g)\in Aut_{\C}({\cal S}) \,\}
$$
denote group of automorphisms of $V$ which preserve $W$ and act on
$Gr^W$ by infinitesimal isometries.  Then, by functoriality, $(F,W)$ 
determines an induced mixed Hodge structure on 
$\gg = Lie(\Gc)$ via the bigrading:
\setcounter{equation}{14}
\begin{equation}
	\gg^{r,s} := gl(V)^{r,s} \cap Lie(\Gc)		
\end{equation}
{\it N.b. By virtue of equation $(2.10)$, 
	  $\overline{\gg^{r,s}} = \gg^{s,r} \mod \oplus_{a<s,b<r}\,\gg^{a,b}$.
\hph{N.b.} Also, $r+s>0\implies\gg^{r,s} = 0$.
}
\end{myremark}

\par To finish our review of Deligne's theory of mixed Hodge structures,
we shall recall the basic properties of the weight filtration $W$ and the
definition of the relative weight filtration 
$$
	\rel W = \rel W(N,W)
$$

\setcounter{definition}{15}				
\begin{definition} Let $W$ be an increasing filtration of a finite dimensional
complex vector space $V$.  Then, a semi-simple endomorphism $Y\in gl(V)$ is 
said to grade $W$ provided that, for each index $k$
$$
	W_k = W_{k-1}\oplus E_k(Y)
$$
(i.e. $W_k$ is the direct sum of $W_{k-1}$ and the $k$-eigenspace of $Y$).
\end{definition}

\setcounter{example}{16}				
\begin{example} Let $(F,W)$ be a mixed Hodge structure.  Then, by definition,
the semi-simple endomorphism $Y$ defined by the rule
$$
	Y(v) = k v \iff v\in\bigoplus_{p+q=k}\, I^{p,q}
$$
is a grading of $W$.
\end{example}

\par To describe the structure of the set of all gradings of a fixed filtration
$W$, let $Lie_{-1}$ denote the nilpotent ideal of $gl(V)$ defined by the rule:
\setcounter{equation}{17}
\begin{equation}
	\a \in Lie_{-1} \iff \a:W_k\to W_{k-1}\hph{a}\forall\, k  
\end{equation}

\setcounter{theorem}{18}			      
\begin{theorem} The unipotent Lie group $\exp(Lie_{-1})$ acts simply 
transitively upon the set of all gradings $Y$ of a fixed, increasing 
filtration $W$.
\end{theorem}
\begin{proof} See \cite{CKS}.
\end{proof}

\par As a prelude to our discussion of the relative weight filtration, let us
consider first a more classical object, namely the monodromy weight
filtration:

\setcounter{theorem}{19}			     
\begin{theorem} Let $V$ be a finite dimensional vector space and $N$ be a 
nilpotent endomorphism of $V$.  Then, there exists a unique 
{\rm monodromy weight filtration} 
$$
	0\subset   W(N)_{-k} \subseteq W(N)_{1-k} \subseteq \cdots
	 \subseteq W(N)_{k-1} \subseteq W(N)_k = V
$$
of $V$ with the following two properties:
\begin{description}
\item{$\bullet$} $N:W(N)_j\to W(N)_{j-2}$ for each index $j$.
\item{$\bullet$} The induced maps 
        	 $N^j:Gr^{W(N)}_{j}\to Gr^{W(N)}_{-j}$ are isomorphisms.
\end{description}
\end{theorem}

\setcounter{example}{20}				
\begin{example} Let $\rho$ be a finite dimensional representation of 
$sl_2(\C)$ and
$$
	N_{\pm} = \rho(n_{\pm}),\qquad Y = \rho(y)
$$
denote the images of the standard generators $(n_-,y,n_+)$ of $sl_2(\C)$.
Then, by virtue of the semi-simplicity of $sl_2(\C)$ and the commutator
relations
$$
	[Y,N_{\pm}] = \pm 2N_{\pm}, \qquad [N_+,N_-] = Y
$$
it follows that: 
$$
	W(N_-)_k = \bigoplus_{j\leq k}\, E_j(Y)
$$
\end{example}

\setcounter{definition}{21}				
\begin{definition} Given an increasing filtration $W$ of a finite dimensional
vector space $V$ and an integer $\ell\in\Z$ the corresponding shifted
object $W[\ell]$ is the increasing filtration of $V$ defined by the rule:
$$
		W[\ell]_j = W_{j+\ell}
$$
\end{definition}

\setcounter{theorem}{22}			     
\begin{theorem} Let $W$ be an increasing filtration of a finite dimensional
vector space $V$.  Then, given a nilpotent endomorphism $N:V\to V$ which 
preserves $W$, there exists at most one increasing filtration
$$
	\rel W = \rel W(N,W)
$$
with the following two properties:
\begin{description}
\item{$\bullet$} \hph{aa} For each index $j$, $N:\rel W_j\to\rel W_{j-2}$
\item{$\bullet$} \hph{aa} For each index $k$, $\rel W$ induces on $Gr^W_k$ the 
		 corresponding shifted monodromy weight filtration
		 $$                  
			W(N:Gr^W_k\to Gr^W_k)[-k]
		 $$
\end{description}
\end{theorem}

\section{Classifying Spaces}
\label{sec:3}


\par In this section, we construct classifying spaces of graded-polarized
mixed Hodge structures $\M$ which generalize Griffiths classifying 
spaces $\D$ of pure, polarized Hodge structures.  In particular, we show
that any variation of graded-polarized mixed Hodge structure $\V\to S$
admits a reformulation in terms of its monodromy representation
$$
	\rho:\pi_1(S,s_0)\to Aut(\V_{s_0})
$$
and its period map
$$
	\phi:S\to\M/\Gamma,\qquad \Gamma = \text{Image}(\rho)
$$

\begin{myremark} Classifying spaces of graded-polarized mixed Hodge structures
have been studied before, notably in \cite{Usui} and the unpublished work
\cite{Kaplan}.  For the most part, the presentation given here follows 
\cite{Kaplan}.
\end{myremark}

\par To establish notation, let us first review Griffiths construction:
Let $V$ be a finite dimensional complex vector space endowed with a rational 
structure $V_{\Q}$ and a non-degenerate bilinear form
$$
	Q:V_{\Q}\otimes V_{\Q}\to\Q
$$
of parity $(-1)^k$.  Then, given any partition of $\text{dim}\,V$ into a sum 
of non-negative integers $\{h^{p,k-p}\}$, one can from the corresponding
classifying space
$$
	\D = \D(V,Q,h^{p,k-p})
$$
consisting of all pure Hodge structure of weight $k$ on $V$ which are
polarized by $Q$ and satisfy
$$
	\text{dim}\,H^{p,k-p} = h^{p,k-p}
$$
 
\par A priori, the classifying space $\D$ is just a set.  To endow it with 
the structure of a complex manifold, one may proceed as follows:  Let
$$
	f^p = \sum_{r\geq p}\, h^{r,k-r}
$$
and $\check\F$ denote the flag variety consisting of all decreasing filtrations
$F$ of $V$ such that
$$
	\text{dim}\,F^p = f^p
$$
Because $\check\F$ is an smooth algebraic manifold, the subset 
$\check\D\subseteq\check\F$ consisting of the those filtrations 
$F\in\check\F$ which satisfy the first Riemann bilinear relation
$$
	Q(F^p,F^{k-p+1}) = 0
$$
is also algebraic.  

\par Now, as can be easily checked via elementary linear algebra, the complex 
Lie group $\Gc=Aut_{\C}(Q)$ acts transitively on $\check\D$.  Consequently,
$\check\D$ is in fact a smooth subvariety of $\check\F$.  

\par To prove that $\D$ is an open subset of $\check\D$, observe that because
$\Gc$ acts transitively on $\check\D$, the map 
$$
	g\in\Gc \mapsto g.F\in\check\D
$$
defines a holomorphic surjection from a neighborhood of $1\in\Gc$ onto a
neighborhood of $F\in\check\D$.  Consequently, by virtue of the following
lemma, there exists an open subset of $\check\D$ about each point 
$F\in\D$ which is entirely contained in $\D$:

\setcounter{equation}{1}					
\begin{lemma} The Lie group $G_{\R} = Aut_{\R}(Q)$ acts transitively on $\D$.
Moreover, given point $F\in\D$, there exists a neighborhood $O$ of $1\in\Gc$ 
such that
\begin{equation}						
	g_{\C}\in O \implies g_{\C}.F \in\D			
\end{equation}
\end{lemma}
\begin{proof} The proof that $G_{\R}$ acts transitively on $\D$ is an exercise
in elementary linear algebra which shall be left to the reader. 

\par To verify equation $(3.2)$, observe that (in the notation of \S 2)
\begin{equation}
	Lie(\Gc) = \bigoplus_p\, \gg^{p,-p}			
\end{equation}
while the Lie algebra of the isotopy group $\Gc^F$ is given by the formula:
\begin{equation}
	Lie(\Gc^F) = \bigoplus_{p\geq 0}\, \gg^{p,-p}		
\end{equation}
Moreover, the subalgebra
$$
	Lie(G_{\R}^F) = Lie(\Gc^F)\cap Lie(G_{\R})\subset Lie(\Gc^F)
$$
consists of exact those elements $\a\in\gg^{0,0}$ which are self-conjugate.
Consequently, 
$$
	{\cal C} = \left(\bigoplus_{p>0}\,\gg^{p,-p}\right)
		 \bigoplus \sqrt{-1}\, Lie(G_{\R}^F)
$$
is a vector space complement to $Lie(G_{\R})$ in $Lie(\Gc)$, i.e.
\begin{equation}
	Lie(\Gc) = Lie(G_{\R}) \oplus {\cal C}			
\end{equation}

\par By virtue of this vector space decomposition, there exists
a neighborhood $U_0$ of zero in $Lie(\Gc)$ such that every element
$$
	g_{\C}\in \exp(U_0)
$$
may be uniquely decomposed into a product
$$
	g_{\C} = g_{\R}g_{\C}^F
$$
of an element $g_{\R}\in G_{\R}$ and an element $g_{\C}^F\in\exp({\cal C})$.
In particular,
$$
	g_{\C}\in\exp(U_0)\implies
	g_{\C}.F = g_{\R}g_{\C}^F.F = g_{\R}.F\in\D
$$ 
\end{proof}

\begin{myremark} A smooth map $F:S\to\D$ is holomorphic provided that relative
to any choice of local holomorphic coordinates $(s_1,\dots,s_n)$ on $S$, one
has
$$
	\frac{\pd F^p}{\pd\bar s_j} \subseteq F^p(s_1,\dots,s_n)
$$
A holomorphic map $F:S\to\D$ is said to be horizontal provided
$$
	\frac{\pd F^p}{\pd s_j} \subseteq F^{p-1}(s_1,\dots,s_n)
$$
\end{myremark}

\par The relationship between variations of pure, polarized Hodge structure
(VHS) and the classifying spaces $\D$ is as follows: Let $\V\to S$ be a 
variation of pure, polarized Hodge structure.  Then, choice of base point 
$s_0\in S$ determines a monodromy representation
\setcounter{equation}{5}
\begin{equation}
	\rho:\pi_1(S,s_0)\to\Gamma				
\end{equation}
and a (locally-liftable) holomorphic, horizontal map
\begin{equation}
	\phi:S\to\D/\Gamma					
\end{equation}
via parallel translation of the data of $\V$ to the reference fiber
$V = V_{s_0}$.  

\par Conversely, given the monodromy representation $(3.6)$ and the period map
$(3.7)$, it is possible to reconstruct the original variation $\V\to S$ 
(up to isomorphism) by reversing the preceding construction.

\par To reformulate the notion of a variation of graded-polarized mixed Hodge 
structure $\V\to S$ in terms of the monodromy representation $\rho$ of $\V$ 
and a suitable period map $\phi:S\to\M/\Gamma$, one must first construct a 
suitable classifying space $\M$ of graded-polarized mixed Hodge structures.  

\par To this end, let $V$ be a complex vector space endowed with a choice of
rational structure $V_{\Q}$ and a choice of weight filtration $W$ (also defined
over $\Q$).  Then, given a collection of rational, non-degenerate bilinear 
forms
$$
	\calS_k:Gr^W_k\otimes Gr^W_k\to\C
$$
of alternating parity $(-1)^k$ and a partition of $\text{dim}\, V$ into 
suitable sum of non-negative integers $\{h^{p,q}\}$, one can form the
classifying space $\M$ consisting of all mixed Hodge structures $(F,W)$ 
which are graded-polarized by $\cal S$ and satisfy the dimensionality 
condition
$$
	\text{dim}\, I^{p,q}_{(F,W)} = h^{p,q}
$$

\par In analogy with the pure case, in order to prove that $\M$ is a complex 
manifold one starts with the flag variety $\check\F$ consisting of all 
decreasing filtrations $F$ of $V$ such that
$$
	\text{dim}\, F^p = f^p,\qquad f^p = \sum_{r\geq p,s}\, h^{r,s}
$$
Next, one defines $\check\F(W)$ to be the submanifold of $\check\F$ consisting 
of those filtrations which have the following additional property:
$$
    \text{dim}\, F^p Gr^W_k = f^p_k,\qquad f^p_k = \sum_{r\geq p}\, h^{r,k-r}
$$
To prove that $\check\F(W)$ is a smooth submanifold of $\check F$, one simply 
checks that the complex Lie group
$$
	GL(V)^W = \{\, g\in GL(V) \mid g:W_k\to W_k\quad\forall k\,\}
$$
acts transitively on $\check\F(W)$.

\par To proceed further, one introduces the \lq\lq compact dual\rq\rq{} 
$\check\M$ consisting of all filtrations $F\in\check\F(W)$ which satisfy the 
first Riemann bilinear relation
$$
	\calS_k(F^p Gr^W_k,F^{k-p+1} Gr^W_k) = 0
$$
for each index $k$.  As in the pure case, in order to show that $\check\M$
is a smooth submanifold of $\check\F(W)$, one proves that a suitable
Lie group acts transitively on $\check\M$:

\setcounter{lemma}{7}						
\begin{lemma} The complex Lie group 
$$
	\Gc = \{\, g\in GL(V)^W \mid Gr(g)\in Aut_{\C}(\calS) \,\}
$$
acts transitively on $\check\M$.	
\end{lemma}
\begin{proof} As in the pure case, Lemma $(3.8)$ can be checked by brute
force using only elementary linear algebra.  Alternatively, one can 
proceed as follows:  Let $\Y(W)$ denote the set of all gradings of $W$,
$$
	\check\D = \bigoplus_k\,\check\D(Gr^W_k,\calS_k,h^{p,k-p})
$$
and $\check\X$ denote the product space $\check\D\times\Y(W)$.  Next,
let 
$$
	\pi:\check\X\to\check\M
$$
denote the natural projection map which sends a point 
$(\oplus_k\,F_k,Y)$ in $\check\X$ to the filtration $F\in\check\M$
determined by the given filtration $\oplus_k\,F^{\cdot}_k$ of $Gr^W$ 
and the induced isomorphism $Y:Gr^W\to V$, i.e. 
$$
	F^p = \bigoplus_k\, Y(F^p_k)
$$

\par Now, as may be easily checked, the map $\pi:\check\X\to\check\M$
is both surjective and $\Gc$ equivariant.  Therefore,  in order to prove 
that $\Gc$ acts transitively on $\check\M$, it will suffice to prove 
that $\Gc$ acts transitively on $\check\X$.  To prove the latter
assertion, observe that given a point $x = (\oplus_k\,F_k,Y)\in\check\X$, 
the corresponding isotopy group $\Gc^Y$ acts transitively on $\check\D$ 
while fixing the given grading $Y$.  On the other hand, the unipotent Lie 
group $\exp(Lie_{-1})$ discussed in \S 2 acts simply transitively on 
$\Y(W)$ while leaving $\D$ pointwise fixed. 
\end{proof}

\par In analogy with the pure case, the proof of the fact that $\M$ is an
open subset of $\check\M$ follows directly from that fact that $\Gc$ acts
transitively on $\check\M$ together with the following lemma:

\setcounter{lemma}{8}						
\begin{lemma} The Lie group
$$
	G = \{\, g\in GL(V)^W \mid Gr(g) \in Aut_{\R}(\calS) \,\}
$$
acts transitively on $\M$.  Moreover, given an element $F\in\M$, there exists 
a neighborhood $U$ of $1\in\Gc$ such that
\setcounter{equation}{9}
\begin{equation}
	g_{\C}\in U \implies g_{\C}.F\in\M			
\end{equation}
\end{lemma}
\begin{proof} The proof of the fact that $G$ acts transitively on $\M$
follows mutatis mutandis from the proof of Lemma $(3.8)$.

\par To verify equation $(3.10)$, pick a point $F\in\M$ and let $Y=Y_{(F,W)}$ 
be the grading defined in Example $(2.17)$.  Then, as may be verified by 
direct computation:
$$
	Lie(\Gc) = Lie(\Gc^Y)\oplus Lie_{-1}
$$
Moreover, the subspace
$$
	{\cal C} = \left(\bigoplus_{p>0}\,\gg^{p,-p}\right)\bigoplus
		   \sqrt{-1}\,\left(Lie(\Gc^F)\cap Lie(G^Y)\right)
$$
is a vector space complement to $Lie(G^Y)$ in $Lie(\Gc^Y)$.  Consequently,
over a sufficiently small neighborhood $U_0$ of zero in $Lie(\Gc)$, every
element
$$
	g_{\C}\in\exp(U_0)
$$
will admit a unique decomposition 
$$
	g_{\C} = g_{-1}g^Yg_{\C}^F
$$
with $g_{-1}\in\exp(Lie_{-1})$, $g^Y\in\exp(Lie(G^Y))$ and $g_{\C}^F$ in 
$\exp({\cal C})\subset\Gc^F$.  In particular,
$$
	g_{\C}\in\exp(U_0) \implies 
	g_{\C}.F = g_{-1}g^Yg_{\C}^F.F = g_{-1}g^Y.F\in\M
$$
\end{proof}

\begin{myremark} $\exp(Lie_{-1})$ is a subgroup of $G$.
\end{myremark}

\par As in the pure case, the relationship between variations of 
graded-polarized mixed Hodge structures (VGPMHS) and the corresponding
classifying spaces of graded-polarized mixed Hodge structures is as follows:
Let $\V\to S$ be a VGPMHS.  Then, choice of a base point $s_0\in S$ 
determines a monodromy representation
\setcounter{equation}{10}
\begin{equation}
	\rho:\pi_1(S,s_0)\to\Gamma			\tag{3.11}
\end{equation}
and a (locally-liftable) holomorphic, horizontal map
\begin{equation}
	\phi:S\to\M/\Gamma				\tag{3.12}
\end{equation}
via parallel translation of the data of $\V$ to the reference fiber 
$V=\V_{s_0}$.

\par Conversely, given the monodromy representation $(3.11)$ and the
period map $(3.12)$, it is possible to reconstruct the original 
variation $\V$ (up to isomorphism) by reversing the preceding 
construction.

\par To close this section, we shall now study the relationship between
$Lie(\Gc)$ and $\M$ in a bit more detail:

\setcounter{theorem}{12}				
\begin{theorem} At each point $F\in\M$, the map
\setcounter{equation}{13}
\begin{equation}
        u\in q_F\mapsto \exp(u).F\in\check\M 		
\end{equation}
restricts to a biholomorphism from a neighborhood of zero in the nilpotent
subalgebra
\begin{equation}
        q_F = \bigoplus_{r<0,r+s\leq 0}\, \gg^{r,s}\subseteq Lie(\Gc) 
\end{equation}
to a neighborhood of $F\in\M$.
\end{theorem}
\begin{proof} Because $\Gc$ acts transitively on $\check\M$ and $\M$
is an open subset of $\check\M$, it suffices to check that $q_F$ is
a vector space complement of $Lie(\Gc^F)$ in $Lie(\Gc)$. To verify
this last assertion, observe that
$$
	Lie(\Gc^F) = \bigoplus_{r\geq 0,r+s\leq 0}\,\gg^{r,s}
$$
and hence
$$
	Lie(\Gc^F)\oplus q_F = \bigoplus_{r+s\leq 0}\,\gg^{r,s} 
			     = Lie(\Gc)
$$
\end{proof}

\setcounter{corollary}{15}				
\begin{corollary} The map $q_F\to T_F(\M)$ which sends the endomorphism
$u\in q_F$ to the derivation
$$
        u(\zeta) = \frac{d}{dt}\zeta(\exp(t u).F)|_{t=0}
$$
is a $\C$-linear isomorphism.  In particular, the derivative $\Phi_{\ast}$
of a holomorphic, horizontal map $\Phi:S\to\M$ takes values in the 
horizontal subbundle
$$
	T^{horiz}_F(\M) = \bigoplus_{k\leq 1}\, \gg^{-1,k}_{(F,W)},
	\qquad F\in\M
$$
\end{corollary}

\begin{myremark} The methods developed above may also be used to show that
the real Lie group
\setcounter{equation}{16}
\begin{equation}
	G_{\R} = \{\, g\in GL(V_{\R})^W \mid Gr(g)\in Aut_{\R}(\calS)\,\}
\end{equation}
acts transitively on the $C^{\infty}$ submanifold $\M_{\R}\subseteq\M$ 
consisting of those filtrations $F\in\M$ for which the corresponding 
mixed Hodge structure $(F,W)$ is split over $\R$, i.e.
$$
	\overline{I^{p,q}_{(F,W)}} = I^{q,p}_{(F,W)}
$$
Moreover, as may be easily checked by direct computation, 
\begin{equation}
	g\in G_{\R},\,F\in\M \implies I^{p,q}_{(g.F,W)} = g.I^{p,q}_{(F,W)}
\end{equation}
In particular, the action of $G_{\R}$ on $\M$ preserves the submanifold 
$\M_{\R}$.
\end{myremark}

\section{Deligne-Hodge bundles}
\label{sec:4}
\setcounter{equation}{0}
\setcounter{definition}{0}
\setcounter{theorem}{0}
\setcounter{lemma}{0}

\par Let $V=V_{\Q}\otimes\C$ be a finite dimensional complex vector space 
which is defined over $\Q$.  Then, as discussed in \S 2, each choice of a 
mixed Hodge structure $(F,W)$ on $V$ determines a unique, functorial 
decomposition
\begin{equation}
	V = \bigoplus_{p,q}\, I^{p,q}_{(F,W)}			
\end{equation}
with the following three properties:
\begin{description}
\item{(1)} $F^p = \oplus_{a\geq p,b}\, I^{a,b}$. 
\item{(2)} $W_k = \oplus_{a+b\leq k}\, I^{a,b}$.
\item{(3)} $\bar I^{p,q} = I^{q,p} \mod \Lambda^{q-1,p-1}(V)$.
\end{description}
Consequently, each classifying space of graded-polarized mixed Hodge $\M$
modeled on $V$ supports a natural decomposition
\begin{equation}
	E = \bigoplus_{p,q}\, \I^{p,q}				
\end{equation}
of the corresponding trivial bundle $E = V\times\M$ into a sum of $C^{\infty}$
subbundles 
\begin{equation}
	\I^{p,q}_F = I^{p,q}_{(F,W)}				
\end{equation}

\par To understand the first order behavior of the decomposition $(4.2)$
relative to the flat connection
$$
	\nabla:\E^0(E)\to \E^1(E)
$$
defined by exterior differentiation, recall that as discussed at the 
end of \S 2, the Lie group
\begin{equation}
	G_{\R} = \{\, g\in GL(V_{\R})^W \mid Gr(g)\in Aut_{\R}(\cal S)\,\}
\end{equation}
act transitively upon the set of \lq\lq real points\rq\rq{} 
$\M_{\R}\subseteq\M$.  Moreover, by virtue of equation $(3.18)$:
$$
	g\in G_{\R}, F\in\M \implies \I^{p,q}_{g.F} = g.I^{p,q}_{(F,W)}
$$
Consequently, it is relatively easy to understand the behavior of the
decomposition $(4.2)$ along
$\M_{\R}$.

\par Unfortunately however, the group $G_{\R}$ does not (in general) act
transitively upon the entire classifying space $\M$.  Therefore, to apply 
the methods of the preceding paragraph to study the local behavior of the 
decomposition $(4.4)$ near a given point $F\in\M$, we must construct
a $C^{\infty}$ decomposition of each element $g_{\C}\in\Gc$ into a product
\begin{equation}
	g_{\C} = g_{\R} \tilde g g_{\C}^F				
\end{equation}
with the following two properties
\begin{description}
\item{$\bullet$} $g_{\R}\in G_{\R}$, $g\in G$, $g_{\C}^F\in\Gc^F$.
\item{$\bullet$} $I^{p,q}_{(\tilde g.F,W)} = \tilde g.I^{p,q}_{(F,W)}$.
\end{description}
In fact, provided that we are only interested in the local behavior of
the Deligne--Hodge bundles, it will suffice to construct $(4.5)$ over
a neighborhood of $1\in\Gc$.

\par Now, as discussed in \cite{Kaplan}, there exists a large class of 
$C^{\infty}$ decompositions of the form $(4.5)$ which are in some
sense \lq\lq natural\rq\rq{}.  However, for the task at hand, the
decomposition determined by the following theorem appears to be
the most suitable:

\setcounter{theorem}{5}						
\begin{theorem} Let $F$ be a point of $\M$ and $\lam_{(F,W)}$ be the 
nilpotent subalgebra of $Lie_{-1}$ defined by the rule
$$
	\lam_{(F,W)} = \bigoplus_{r,s<0}\,\gg^{r,s}
$$
Then,
$$
	g\in\exp(\lam_{(F,W)})\implies \I^{p,q}_{g.F} = g.\I^{p,q}_F
$$
Moreover, there exists a natural $\R$-vector subspace $\Phi_F\subset Lie(\Gc)$
such that
\setcounter{equation}{6}
\begin{equation}
	Lie(\Gc) = Lie(G_{\R})\oplus
		   \sqrt{-1}\left(\lam_{(F,W)}\cap Lie(G_{\R})\right)\oplus
		   \Phi_F
\end{equation}
\end{theorem}

\setcounter{corollary}{7}				
\begin{corollary} Let $F$ be an element of $\M$.  Then, there exists a 
neighborhood $\exp(U_0)$ about $1\in\Gc$ such that each element 
$g_{\C}$ in $\exp(U_0)$ admits a unique, $C^{\infty}$ decomposition
\setcounter{equation}{9}
\begin{equation}
	g_{\C}= g_{\R}e^{\l}e^{\phi}	
\end{equation}
such that
\begin{description}
\item{$\bullet$} $g_{\R}$ is an element of $G_{\R}$. 
\item{$\bullet$} $e^{\l}$ is an element of 
		 $\exp(\sqrt{-1}\,\lam_{(F,W)}\cap Lie(G_{\R}))$.
\item{$\bullet$} $\exp(\phi)$ is an element of $\exp(\Phi_F)$.
\end{description}
Moreover, along the subalgebra $q_F$, 
$$
	\phi(u) = -\pip(\bar u) 
	+ \text{(higher order terms in $u$ and $\bar u$)}
$$
\end{corollary}
 
\par The proof of Theorem $(4.6)$ and Corollary $(4.8)$ will occupy the
remainder of this section.  In essence however, the proofs of these two
results boil down to a series of relatively straightforward calculations.

\setcounter{lemma}{9}						
\begin{lemma} Let $F$ be a point of $\M$.  Then, the corresponding subgroup
$$
	\exp(\lam_{(F,W)})\subseteq\exp(Lie_{-1})
$$
is closed under conjugation.
\end{lemma}
\begin{proof}  As discussed in \S 2,
$\overline\gg^{r,s} = \gg^{s,r} \mod \oplus_{a<s,b<r}\,\gg^{a,b}$.
Consequently, 
\begin{eqnarray*}
	\overline{\lam_{(F,W)}} &=& \overline{\bigoplus_{r,s<0}\,\gg^{r,s}}
		    = \bigoplus_{r,s<0}\,\overline\gg^{r,s} 	
		    = \bigoplus_{r,s<0}\,
		       \left(\gg^{s,r} \mod \bigoplus_{a<s,b<r}\,\gg^{a,b}
		       \right) \\
		   &=& \bigoplus_{r,s<0}\,\gg^{s,r} 
		    = \lam_{(F,W)}
\end{eqnarray*}
\end{proof}

\setcounter{lemma}{10}						
\begin{lemma} Let $F$ be a point of $\M$.  Then, 
$$
	g\in\exp(\lam_{(F,W)})\implies \I^{p,q}_{g.F} = g.\I^{p,q}_F
$$
\end{lemma}
\begin{proof} 
Let $\{J^{p,q}\}$ denote the bigrading of $(g.F,W)$ 
defined by the rule 
$$
	J^{p,q} = g.I^{p,q}_{(F,W)}
$$
and note that, by construction, each element $h\in\exp(\lam_{(F,W)})$ preserves
the condition 
$$
	v\in I^{q,p}_{(F,W)} \mod \Lambda^{q-1,p-1}(V)
$$
In particular, because the subgroup $\exp(\lam_{(F,W)})$ is closed under
conjugation:
\begin{eqnarray*}
	\bar J^{p,q} &=& \bar g.\bar I^{p,q}
	 = \bar g.\left(I^{q,p} \mod \Lambda^{q-1,p-1}(V)\right) \\
	&=& g\left(g^{-1}\bar g\right).
	    \left(I^{q,p} \mod \Lambda^{q-1,p-1}(V)\right)   \\
	&=& g.\left(I^{q,p} \mod \Lambda^{q-1,p-1}(V)\right) \\ 
	&=& J^{q,p} \mod \bigoplus_{a<q,b<p}\, J^{a,b}
\end{eqnarray*}
Thus, by uniqueness, $J^{p,q} = I^{p,q}_{(g.F,W)}$.
\end{proof}

\par Armed with these two preliminary lemmata, we are now ready to complete
the proofs of Theorem $(4.6)$ and Corollary $(4.8)$:

\begin{proof}[Theorem (4.6)] Observe first that each point $F\in\M$ determines
a vector space decomposition of $Lie(\Gc)$ into direct sum of subalgebras
\setcounter{equation}{11}
\begin{equation}
\begin{array}{ll} 
   \np = \bigoplus_{r\geq 0,\,s<0}\,\gg^{r,s}, & \qquad\nz = \gg^{0,0} \\
   \nn = \bigoplus_{s\geq 0,\,r<0}\,\gg^{r,s}, & \qquad\lam = \lam_{(F,W)}
\end{array}							
\end{equation}
with the following properties:
\begin{equation}
\begin{array}{ll}
     q_F = \nn\oplus\lam_F, & \qquad Lie(\Gc^F) = \np\oplus\nz \\
     \bar\np\subseteq\nn\oplus\lam, & \qquad \bar\nz\subseteq\nz\oplus\lam \\
     \bar\nn\subseteq\np\oplus\lam, & \qquad \bar\lam = \lam        
\end{array}
\end{equation}

\par Next, let $\pip$, $\piz$, $\pin$ and $\pil$ denote projection from
$Lie(\Gc)$ to the corresponding subalgebras $\np$, $\nz$, $\nn$ and $\lam$
listed in equation $(4.12)$, and define
\begin{equation}
  \Phi_F = \np\oplus\{\, x\in\nz\mid \piz(\bar x) = -\piz(x)\,\} 
\end{equation}

\par Finally, observe that since equation $(4.7)$ is a linear condition,
it will suffice to check its validity on each of the subalgebras appearing
in equation $(4.12)$.  Direct computation shows that:
\begin{eqnarray*}
&x &\in\np \implies x = [0]\oplus[0]\oplus[x] \\ 
&x &\in\nz \implies x = [\Re(x)]\oplus[\pil(\Im(x))]\oplus[\piz(\Im(x))] \\
&x &\in\nn \implies x = [\Re(2x-\pil(\bar x))]
                        \oplus [-\Im(\pil(\bar x))]    
                        \oplus [-\pip(\bar x)]  \\
&x &\in\lam \implies x = [\Re(x)]
                           \oplus[\Im(x)]\oplus[0]
\end{eqnarray*}
where as usual,
$$
	\a\in Lie(\Gc) \implies\left\{
	\begin{array}{l} 
	\Re(\a) = \half(\a+\bar\a) \\
	\Im(\a) = \half(\a-\bar\a)
	\end{array}\right.
$$
For example, $x\in\nn\implies \bar x = \pil(\bar x) + \pip(\bar x)$
and hence
\begin{eqnarray*}
\lefteqn{\Re(2x - \pil(\bar x)) - \Im(\pil(\bar x) - \pip(\bar x))
         \hph{aaaaaaaaaaaaaaaaaaaaaaaaaaaaaaaaaaa} }  \\  
  &=& x + \bar x - \half(\pil(\bar x) + \overline{\pil(\bar x)})  
   -\half(\pil(\bar x) - \overline{\pil(\bar x)}) - \pip(\bar x) \\
  &=& x + \bar x - \pil(\bar x) - \pip(\bar x) 	= x
\end{eqnarray*}
\end{proof}

\begin{proof}[Corollary (4.8)]  Let $u$ be an element of $U_0\cap q_F$
$$
	e^u = g_{\R}(u)e^{\l(u)}e^{\phi(u)}
$$
denote the decomposition of the element $e^u\in\exp(U_0)$ defined by 
equation $(4.9)$, and define $\g(u)\in Lie(G_{\R})$ by the rule
$$
	g_{\R}(u) = e^{\g(u)}
$$
Then, applying the Campbell--Baker--Hausdorff formula, one finds that
\begin{equation}
	e^u = e^{\g(u)}e^{\l(u)}e^{\phi(u)} 
	    = e^{\g(u)+\l(u)+\phi(u) + \text{(higher order brackets)}}
\end{equation}

\par To use equation $(4.15)$ to determine the first order behavior of 
$\phi(u)$, one simply inserts the first order Taylor series expansions
$$
	\g(u) = \g_1(u) + O^2(u),\quad \l(u) = \l_1(u) + O^2(u), \quad
	\phi(u) = \phi_1(u) + O^2(u)
$$
into expression $(4.15)$, thereby obtaining the equation
\begin{equation}
	e^u = e^{\g_1(u) + \l_1(u) + \phi_1(u) + O^2(u)}	
\end{equation}
Comparing the linear terms on each side of equation $(4.16)$, it therefore
follows that
\begin{equation}
	u = \g_1(u) + \l_1(u) + \phi_1(u)			
\end{equation}
Applying our previous formulae to equation $(4.17)$, and remembering that $u$ 
is an element of $q_F$, we obtain the desired result.  Namely:
$$
	\phi_1(u) = -\pip(\bar u)
$$
\end{proof}

\begin{myremark} Throughout this paper I shall use the symbol $O^2(x)$ to 
denote an error term of order 2 depending (in principle) upon both $x$ and 
$\bar x$.
\end{myremark}

\section{Higgs Fields}
\label{sec:5}

\par In this section, we prove the main theorem of this paper, namely:

\setcounter{theorem}{0}					      
\begin{theorem} Let $V\to S$ be a variation of graded-polarized mixed Hodge 
structure, and $\{\U^p\}$ denote the collection of $C^{\infty}$ subbundles 
of $\V$ defined by the rule:
\setcounter{equation}{1}
\begin{equation}
	\U^p_s = \bigoplus_q\, I^{p,q}_{(\F_s,\W_s)}	
\end{equation}
Then, relative to the Gauss--Manin connection $\nabla$, the direct sum 
decomposition
$$
	\V = \bigoplus_p\,\U^p				
$$
defines a (unpolarized) complex variation of Hodge structure.
\end{theorem}

\par In particular, by virtue of our discussions in \S 1, the preceding 
result has the following immediate corollary:

\setcounter{corollary}{2}				
\begin{corollary}  Every variation of graded-polarized mixed Hodge
structure $V$ supports a natural Higgs bundle structure $(\V,\bar\pd+\theta)$.
Moreover, because this Higgs bundle structure arises from a complex
variation of Hodge structure, it is automatically a fixed point of the
$\C^*$ action:
$$
		(\V,\bar\pd+\theta) \mapsto (\V,\bar\pd + \l\theta)
$$
\end{corollary}

\par The formal proof of Theorem $(5.1)$ presented below depends upon a 
couple of technical computations [namely: Lemma (5.11)].  The gist of 
the proof however is relative simple, and may be outlined as follows: 
The triple $(\V,\{\U^p\},\nabla)$ defines a complex variation of Hodge 
structure if and only if differentiation induces a map
\setcounter{equation}{3}
\begin{equation}
 \nabla:\E^0(\U^p)\to\E^{0,1}(\U^{p+1})\oplus\E^1(\U^p)\oplus\E^{1,0}(\U^{p-1})
\end{equation}
Thus, in order to prove Theorem $(5.1)$, it will suffice to compute the 
derivative of an arbitrary $C^{\infty}$ local section $\sigma$ of $\U^p$
at a given point $\underline s\in S$.

\par In particular, because the value of $\nabla\sigma$ at $\underline s$ 
is completely determined by the local behavior of $\V$, we may assume that 
our variation is defined over the polydisk
$$
	\Delta^n = \{ (s_1,\dots,s_n)\in\C^n \mid |s_j|<1
		      \hskip 10pt j=1,\dots, n\,\}
$$
via a holomorphic, horizontal map
\begin{equation}
	F(s):\Delta^n\to\M					
\end{equation}
We may also assume that our given point $\underline s\in S$ corresponds to the
point $0=(0,\dots,0)\in\Delta^n$.  

\par In particular, by virtue of Theorem $(3.13)$, there exists a 
neighborhood $O$ about $0\in\Delta^n$ over which the period map $(5.5)$ 
admits a unique representation
\begin{equation}
	F(s) = e^{\Gamma(s)}.F(0)				
\end{equation}
relative to a holomorphic function $\Gamma(s)$ which take values in $q_{F(0)}$
and vanishes at zero.

\par To compute $\nabla\sigma$, observe that in light of Corollary $(4.8)$, 
we may decompose the function $e^{\Gamma(s)}$ into a product of three factors
\begin{equation}
	e^{\Gamma(s)} = g_{\R}(s) e^{\l(s)} e^{\phi(s)}		
\end{equation}
such that
\begin{description}
\item{$\bullet$} $g_{\R}(s)$ takes values in $G_{\R}$.
\item{$\bullet$} $e^{\l(s)}$ takes values in $\exp(\sqrt{-1}\lam_{(F(0),W)}
					 \cap Lie(G_{\R}))$.
\item{$\bullet$} $e^{\phi(s)}$ takes values in $\exp(\Phi_{F(0)})$.
\end{description}
Consequently, the section $\sigma(s)$ may be written in the form
\begin{equation}
	\sigma(s) = g_{\R}(s)e^{\l(s)}.\tilde\sigma(s) 		
\end{equation}
relative to a smooth function $\tilde\sigma(s)$ which takes values in the
fixed vector subspace
\begin{equation}
	\U^p_{F(0)}= \bigoplus_q\,I^{p,q}_{(F(0),W)} 	
\end{equation}		

\par By Leibniz's rule:
\begin{equation}
	\nabla \sigma(s) 
	 	= \left(\nabla g_{\R}(s)e^{\l(s)}\right).\tilde\sigma
		  + \left(g_{\R}(s)e^{\l(s)}\right).\nabla \tilde\sigma
\end{equation}
Thus, to complete the proof of Theorem $(5.1)$, it will suffice to compute 
the value of equation $(5.10)$ at $s=0$.  To this end, we shall employ 
the following lemma:

\setcounter{lemma}{10}					
\begin{lemma} Let $F$ be a point of $\M$ and $u$ be an element of 
$T^{horiz}_F(\M)$.  Then, relative to isomorphism $q_F\cong T_F(\M)$ 
determined by Corollary $(3.16)$,
$$
	\pip(\bar u)\in \gg^{1,-1}\oplus
	     \left(\bigoplus_{k\leq -1}\,\gg^{0,k}\right)
$$
In particular,
$$
	\pip(\bar u):\U^p_{F(0)}\to \U^{p+1}_{F(0)}\oplus\U^p_{F(0)}
$$
\end{lemma}
\begin{proof} Given element $\a\in Lie(\Gc)$ and a point $F\in\M$, let 
$$
	\a = \sum_{r+s\leq 0}\,\a^{r,s},\qquad \a^{r,s}\in\gg^{r,s}
$$
denote the decomposition of $\a$ according to the induced bigrading
$$
	Lie(\Gc) = \bigoplus_{r+s\leq 0}\,\gg^{r,s}
$$
Then, by Corollary $(3.16)$,
$$
	u\in T^{horiz}_F(\M) \implies u =\sum_{k\leq 1}\, u^{-1,k}
$$
Thus,
\setcounter{equation}{11} 
\begin{equation}
	u = u^{-1,1} + u^{-1,0} + \sum_{k\leq -1}\, u^{-1,k}
\end{equation}
In particular, because the rightmost term $\sum_{k\leq -1}\, u^{-1,k}$ of 
equation $(5.12)$ is an element of $\lam_{(F,W)}$ and 
$\lam_{(F,W)} = \overline{\lam_{(F,W)}}$, we have:
$$
	\bar u = \overline{u^{-1,1}} + \overline{u^{-1,0}} \mod \lam_{(F,W)}
$$
To finish the proof, recall that 
$$
	\overline \gg^{r,s} = \gg^{s,r} \mod \bigoplus_{a<s,b<r}\,\gg^{a,b}
$$
Consequently,
\begin{equation}
	\bar u \in 
	   \gg^{1,-1}\oplus\left(\sum_{k\leq -1}\,\gg^{0,k}\right) 
	   \mod \lam_{(F,W)}
\end{equation}
\end{proof}

\begin{proof}[Theorem (5.1)] Let $\nabla = \nabla^{0,1} + \nabla^{1,0}$
denote the decomposition of the Gauss-Manin connection $\nabla$ into its
holomorphic and anti-holomorphic parts.  Then, by virtue of equation $(5.10)$,
Lemma $(5.11)$, and the fact that 
$$
   \Gamma(0) = 0 \implies g_{\R}(0)=1,\quad e^{\l(0)} = 1,\quad
		 \tilde\sigma(0) = \sigma(0)
$$
it will suffice to show that [relative to the isomorphism 
$T_{F(0)}(\M)\cong q_{F(0)}\subset Lie(\Gc)$]:
\begin{equation}
	\begin{array}{ll}
	  \left.\nabla^{0,1}\, g_{\R}(s)e^{\l(s)}\right|_{s=0}
	      & \in \pip\left(\overline{T^{horiz}_{F(0)}(M)}\right)
		    \otimes T^{0,1}_0(\Delta^n)^* \\
	  \left.\nabla^{1,0}\, g_{\R}(s)e^{\l(s)}\right|_{s=0}
	     & \in T^{horiz}_{F(0)}(\M)\otimes T^{1,0}_0(\Delta^n)^* 
	\end{array}					
\end{equation}

\par To this end, let 
\begin{equation}
	\Gamma(s) = \xi_1 s_1 + \cdots + \xi_n s_n + O^2(s)	
\end{equation}
denote the first order Taylor series expansion of $\Gamma(s)$.  Then,
equation $(5.6)$ together with the horizontality of $\V$ imply that 
\begin{equation}
	\xi_j \in T^{horiz}_{F(0)}(\M),\qquad j=1,\dots,n	
\end{equation}

\par Next, observe that by virtue of equation $(5.7)$,
\begin{equation} 
	g_{\R}(s)e^{\l(s)} = e^{\Gamma(s)}e^{-\phi(s)}		
\end{equation}
Moreover, equations $(5.15)$ and $(5.16)$ together with Corollary $(4.8)$,
imply that 
\begin{equation}
	\phi(s) = -\sum_{j=1}^n\,\pip(\bar\xi_j) \bar s_j + O^2(s) 
\end{equation}
Thus, by the Campbell--Baker--Hausdorff formula,
\begin{equation}
 	e^{\Gamma(s)}e^{-\phi(s)} 
	= \exp(\sum_{j=1}^n\,\xi_j s_j 
		+ \sum_{j=1}^n\,\pip(\bar\xi_j) \bar s_j + O^2(s))
\end{equation}
In particular,
\begin{equation}
	\begin{array}{ll}
	 \left.\nabla^{0,1}\, g_{\R}(s)e^{\l(s)}\right|_{s=0} 
	 & = \sum_{j=1}^n\, \pip(\bar\xi_j)\otimes d\bar s_j \\
	 \left.\nabla^{1,0}\,  g_{\R}(s)e^{\l(s)}\right|_{s=0}
	 & = \sum_{j=1}^n\, \xi_j\otimes ds_j 
	\end{array}	
\end{equation}
\end{proof}

\setcounter{corollary}{20}				
\begin{corollary} The Higgs bundle structure $\bar\pd+\theta$ associated to
a variation of graded-polarized mixed Hodge structure $\V\to S$ preserves
the weight filtration $\W$.
\end{corollary}

\par To state the next result, let $(E,\{\U^p\},\nabla)$ be a complex 
variation of Hodge structure.  Then, as discussed in \S 1, the corresponding 
Higgs bundle $(E,\bar\pd+\theta)$ is obtained by using horizontality condition
$$
	\nabla:\E^0(\U^p)\to \E^{0,1}(\U^{p+1})\oplus\E^{0,1}(\U^p)\oplus
			     \E^{1,0}(\U^p)\oplus\E^{1,0}(\U^{p-1})
$$
to decompose $\nabla$ into the sum of a pair of differential operators
$$
	\bar\pd:\E^0(\U^p)\to \E^{0,1}(\U^p),\qquad
	    \pd:\E^0(\U^p)\to \E^{1,0}(\U^p)
$$
and a pair of tensor fields
$$
	\tau   \in Hom(\U^p,\U^{p+1})\otimes\E^{0,1},\qquad
	\theta \in Hom(\U^p,\U^{p-1})\otimes\E^{1,0}
$$

\begin{myremark} During the remainder of this paper, we shall also use the
following notation:  Let $\a$ be an element of $Lie(\Gc)$ and $y$ be a 
grading of $W$.  Then, $\a^Y$ will denote the component of $\a$ which 
is of weight zero relative to the semi-simple endomorphism $ad\,Y$.
\end{myremark}

\setcounter{lemma}{21}					
\begin{lemma} Let $\V$ be a variation of graded-polarized mixed Hodge
structure with Deligne grading 
$$
	\Y(\sigma) = k \sigma \iff 
	\sigma\in\E^0\left(\bigoplus_{p+q = k}\,\I^{p,q}\right)
$$
and Gauss--Manin connection $\nabla = \tau + \bar\pd + \pd + \theta$, 
Then, 
\setcounter{equation}{22}
\begin{equation}
        \tau = (\bar\theta)^{\Y} 
\end{equation}
\end{lemma}

\setcounter{corollary}{23}				
\begin{corollary} Let $\V\to S$ be a VGPMHS.  Then, $\U^p$ is a holomorphic 
subbundle of $\V$ relative to the integrable complex structure 
$\nabla^{0,1}$ iff $\bar\theta^Y$ vanishes on $\U^p$.  In particular, 
assuming the system of Hodge bundles is of the form
\setcounter{equation}{24}
\begin{equation}
    \U^a\oplus\U^{a+1}\oplus\cdots\oplus\U^b 
\end{equation}
it follows that $\U^b$ is a holomorphic subbundle of $\V$ (relative to 
the Gauss--Manin connection $\nabla$).
\end{corollary}

\par The remainder of this section is devoted to some applications of 
Theorem $(5.1)$ to the study of unipotent variations of mixed Hodge
structure.  Applications of Theorem $(5.1)$ to quantum cohomology and
mirror symmetry will be discussed in \S 7 and \S 8.

\setcounter{definition}{25}				
\begin{definition} A variation $\V$ is unipotent if and only if the 
induced variations ${\F}Gr^{\W}$ are constant.
\end{definition}

\setcounter{lemma}{26}					
\begin{lemma} If $\V$ is unipotent then 
\setcounter{equation}{27}
\begin{equation}
        \theta:\W_k\to\W_{k-1} 
\end{equation}
for each index $k$.  In particular, $(\bar\theta)^{\Y} = 0$.
\end{lemma}
\begin{proof} The induced maps $\theta^Y:Gr^{\W}_k\to Gr^{\W}_k$ coincide with
the Higgs field carried by the corresponding variation on $Gr^{\W}_k$.
However, by unipotency, the induced variations on $Gr^W$ are constant.
\end{proof}

\setcounter{corollary}{28}				
\begin{corollary} If $\V$ is a unipotent VGPMHS then 
\begin{description}
\item{(1)}  The complex structures $\bar\pd$ and $\nabla^{0,1}$ coincide.
\item{(2)}  The connection $\bar\pd + \pd = \nabla-\theta$ is flat.
\item{(3)}  The Higgs field $\theta$ is flat relative both to $\nabla$ and 
            $\bar\pd + \pd$.
\end{description}
\end{corollary}

\setcounter{theorem}{29}				
\begin{theorem} A unipotent VGPMHS $\V\to S$ may be recovered from the 
following data:
\begin{description}
\item{(1)}  The flat connection $\nabla$ of $\V$.
\item{(2)}  The associated Higgs field $\theta$.
\item{(3)}  A single fiber $\L_{s_0}$ of $\V$.
\end{description}
\end{theorem}
\begin{proof} Since the weight filtration $\W$ and the bilinear forms
$\{\calS_k\}$ are flat, the key step is to recover the Hodge filtration
$\F$ via the subbundles $\{\U^p\}$.  However, by $(5.29)$, $\U^p$ is parallel
with respect to the flat connection $\nabla-\theta$.
\end{proof}

\par Now, as observed by Deligne \cite{Deligne}, the class of unipotent 
variations described above contains the following special subclass 
which appears to be of great importance in mirror symmetry:

\setcounter{definition}{30}				
\begin{definition} A variation $\V\to S$ is said to be Hodge--Tate provided 
the graded Hodge numbers $h^{p,q}$ vanish unless $p=q$.
\end{definition}

\setcounter{lemma}{31}					
\begin{lemma} Let $\V$ be a variation of Hodge--Tate type.  Then, the 
corresponding Higgs field $\theta$ assumes values $\gg^{-1,-1}$.
\end{lemma}
\begin{proof} In this case, the horizontal subspace 
$$
	T^{horiz}_{F(0)}(\M)\subseteq T_{F(0)}(\M)
$$
appearing in the proof of Theorem $(5.1)$ reduces to $\gg^{-1,-1}_{(F(0),W)}$.
\end{proof}

\setcounter{corollary}{32}				
\begin{corollary} If $\V$ is a Hodge--Tate variation then 
\begin{description}
\item{(1)}  The grading $\Y$ is flat relative to the connection 
            $\bar\pd + \pd$.
\item{(2)}  $\nabla\Y = 2\theta$.
\end{description}
\end{corollary}

\section{Asymptotic Behavior}
\label{sec:6}
\setcounter{equation}{0}

\par In this section we consider the asymptotic behavior of variations 
of graded-polarized mixed Hodge structure $\V\to\Delta^{*n}$ which are 
admissible in the sense of \cite{S--Z} and \cite{Kashiwara}.  
\vskip 10pt

\par To this end, recall that given a flat vector bundle $E\to\Delta^{*n}$ 
with unipotent monodromy, there exists a canonical extension $E^c\to\Delta^n$ 
relative to which the flat connection $\nabla$ of $E$ has at worst a simple 
poles with nilpotent residues along the divisor $D=\Delta^n-\Delta^{*n}$.  

\par More explicitly, given a choice of coordiantes $\Delta^n$ relative to 
which 
\begin{equation}
        D = \{\, p\in\Delta^n \mid s_1(p)\cdots s_n(p) = 0 \,\} 
\end{equation}
one may identify $E^c$ with the locally free sheaf generated by the sections
\begin{equation}
     \sigma^c = e^{\sum_j\, \frac{1}{2\pi i}(\log s_j)N_j}\sigma,\quad
     \sigma \hph{a}\text{a flat, multivalued section of}\hph{a} E
\end{equation}
where $T_j(s):E_s\to E_s$ denotes the action of parallel translation along 
the counterclockwise path
$\g_j(t) = (s_1,\dots,s_{j-1}, e^{2\pi i t} s_j, s_{j+1},\dots,s_n)$. 
and 
$$
	N_j = -\log T_j
$$
In particular, given the choice of coordinates $(6.1)$, we obtain a flat
coonection $\nabla^c$ on $E^c$ via the rule:
$$
        \nabla^c =  \nabla 
        - \frac{1}{2\pi i}\sum_{j=1}^n\,\frac{ds_j}{s_j}\otimes N_j,\qquad 
$$

\begin{myremark} Our sign conventions regarding $N_j$ follow \cite{Deligne}.
\end{myremark}

\par Suppose now that $\V\to\Delta^{*n}$ is a variation of graded-polarized 
mixed Hodge structure with unipotent monodromy.  Then, by virtue $(6.2)$, we 
obtain
\begin{description}
\item{(1)} \hph{aa} A choice of rational structure $V_{\Q}$ on the central 
	   fiber of $\V^c$.
\end{description}
Relative $V_{\Q}$, parallel translation to the central fiber of $E^c$ under 
$\nabla^c$ then defines: 
\begin{description}
\item{(2)} \hph{aa} A rational weight filtration $W$ of $V$.
\item{(3)} \hph{aa} Rational, non-degenerate bilinear forms 
           $\calS_k:Gr^W_k\otimes Gr^W_k\to\C$ of parity $(-1)^k$.
\item{(4)} \hph{aa} An \lq\lq untwisted\rq\rq{} period map 
           $\psi:\Delta^{*n}\to\check\M$ encoding the Hodge filtration 
           $\F$ of $\V$.
\end{description}
\pagebreak

\par Following \cite{S--Z} and \cite{Kashiwara} let us henceforth assume 
that $\V$ is admissible, i.e.
\begin{description}
\item{(5)} \hph{aa} The map $\psi:\Delta^{*n}\to\check\M$ extends to 
	   $\Delta^n$.
\item{(6)} \hph{aa} The data $(1)$--$(3)$ together with the limiting filtration
           \begin{equation}
              F_{\infty} := \lim_{p\to 0}\,\psi(p) 
           \end{equation}
           and the monodromy logarithms $N_1,\dots,N_n$ define an 
           infinitesimal mixed-Hodge module (IMHM) in the sense of 
           \cite{Kashiwara}.
\end{description}
Then, by the work of Deligne and Kashiwara, the monodromy cone
\begin{equation}
        {\cal{C}} = \{\, \sum_{j=1}^n\,a_j N_j \mid a_j>0 \,\}  
\end{equation}
of $\V$ enjoys the following properties:
\begin{description}
\item{$\bullet$} \hph{a} The relative weight filtration $\rel W = \rel W(N,W)$
                 is constant on $\cal{C}$.
\item{$\bullet$} \hph{a} The pair $(F_{\infty},\rel W)$ is a mixed Hodge 
		 structure.
\item{$\bullet$} \hph{a} Each element $N_j$ is a $(-1,-1)$ morphism
                 of $(F_{\infty},\rel W)$.
\end{description}
Consequently, the limiting mixed Hodge structure $(F_{\infty},\rel W)$
of $\V$ defines a canonical decomposition 
\begin{equation}
        Lie(\Gc) = \bigoplus_a\,\wp_a,\qquad 
	\wp_a:U^p_{\infty}\to U^{p+a}_{\infty}
\end{equation}
via the subspaces
\begin{equation}
    U^p_{\infty} = \bigoplus_{q}\, I^{p,q}_{(F_{\infty},\rel W)} 
\end{equation}

\par In particular, since $F_{\infty}^p = F^{p+1}_{\infty}\oplus U^p_{\infty}$,
 the graded, nilpotent Lie algebra
\begin{equation}
        q_{\infty} = \bigoplus_{a<0}\,\wp_a 
\end{equation}
is a vector space complement to $Lie(\Gc^{F_{\infty}})$ in $Lie(\Gc)$.
As a result, there exists a polydisk $\Delta^n_r$ of radius $r$ over 
which 
\begin{equation}
        \psi(s) = e^{\Gamma(s)}.F_{\infty} 
\end{equation}
relative to a unique holomorphic function $\Gamma:\Delta^n_r\to q_{\infty}$
vanishing at zero. 

\begin{myremark} Since the original period map $F(z)$ of $\V$ may be
written as $F(z) = e^{zN}e^{\Gamma(s)}.F_{\infty}$, equation $(6.8)$
may be viewed as a type of \lq\lq group theoretic\rq\rq{} version of
Schmid's Nilpotent Orbit Theorem for admissible VGPMHS.
\end{myremark}

\setcounter{theorem}{8}					
\begin{theorem} The function $\Gamma$ defined by equation $(6.8)$ 
satisfies the differential equation
\setcounter{equation}{9}
\begin{equation}
        e^{-ad\,\Gamma}\,\Omega
        + e^{-\Gamma}\,\pd e^{\Gamma}\in\wp_{-1},\quad
        \Omega = \frac{1}{2\pi i}\sum_{j=1}^n\,\frac{ds_j}{s_j}\otimes N_j
\end{equation}
\end{theorem}
\begin{proof} Let
$$
        U^n_r = \{\, (z_1,\dots,z_n)\in\C^n \mid 
                     (e^{2\pi i z_1},\dots,e^{2\pi i z_n})\in\Delta^{*n}_r \,\}
$$
cover $\Delta^{*n}_r$ via $s_j = e^{2\pi i z_j}$.  Then, the horizontality
of the Hodge filtration
\begin{equation}
   F(z_1,\dots,z_n) = e^{\sum_{j=1}^n z_j N_j}
                      e^{\Gamma(s_1,\dots,s_n)}.F_{\infty}
\end{equation}
implies that 
\begin{equation}
        \left[e^{\sum_{j=1}^n z_j N_j}e^{\Gamma}\right]^{-1} 
        \pd\left[e^{\sum_{j=1}^n z_j N_j}e^{\Gamma}
           \right]\in\wp_{-1}
\end{equation}
which unravels to yield $(6.10)$.
\end{proof}

\setcounter{theorem}{12}			
\begin{theorem} Let $X_{-k}$ denote the component of 
$$
        X = \log\left(e^{\sum_{j=1}^n z_j N_j}e^{\Gamma(s_1,\dots,s_n)}\right)
$$
taking values in $\wp_{-k}$.  Then, the endomorphism valued 1-form
\setcounter{equation}{13}
\begin{equation}
        \pd X_{-1} = \Omega + \pd\Gamma_{-1} 
\end{equation}
defines a Higgs field on the product bundle 
$V\times\Delta^{*n}_r\to\Delta^{*n}_r$.
\end{theorem}
\begin{proof} With a little work, equation $(6.13)$ can be recast as
$$
        \pd e^{X} = e^{X}\,\pd X_{-1}    
$$
Consequently, 
$$
        e^{-X}\frac{\pd^2 }{\pd z_i\pd z_j} e^X 
        = \frac{\pd X_{-1}}{\pd z_i}\frac{\pd X_{-1}}{\pd z_j}
          + \frac{\pd ^2 X_{-1}}{\pd z_i\pd z_j} 
$$
and hence
$$
        \pd X_{-1}\wedge\pd X_{-1} = 0    
$$
by equality of mixed partial derivatives.
\end{proof}

\begin{myremark} As a result of Theorem $(6.13)$,
\begin{equation}
        \tilde\nabla\sigma = d\sigma + \pd X_{-1}(\sigma) 
\end{equation}
defines a flat connection on the product bundle 
$V\times\Delta^{*n}_r\to\Delta^{*n}_r$.
\end{myremark}

\par Comparing Theorem $(6.13)$ with Theorem $(5.30)$, it is natural to ask
when a holomorphic function $\Gamma_{-1}:\Delta^n_r\to\wp_{-1}$ 
defines a solution to $(6.10)$.

\setcounter{theorem}{15}			
\begin{theorem} There exists a bijective correspondence between
solutions to $(6.10)$ and holomorphic functions 
$\Gamma_{-1}:\Delta^n_r\to\wp_{-1}$ which satisfy the Higgs field
condition $\pd X_{-1}\wedge\pd X_{-1}=0$ and the initial value constraint 
$\Gamma_{-1}(0)=0$.
\end{theorem}

\par I shall sketch two proofs of this result.  The first, inspired
by \cite{Deligne}, is to show that $\Gamma$ may be recovered from the 
monodromy action of the flat connection $(6.15)$.  More precisely, the 
monodromy of $(6.15)$ about the loop 
$\g_j(t) = (s_1,\dots,e^{2\pi i t}s_j,\dots,s_n)$, $0\leq t\leq 1$ is given by
\setcounter{equation}{16} 
\begin{equation}
        T_j = e^{-\Gamma}e^{-N_j}e^{\Gamma}    
\end{equation}      

\par The second method of proof, used by Cattani and Peters in their study
of variations of pure Hodge structure, is to write
\begin{equation}
        e^{\Gamma} = 1 + E_{-1} + E_{-2} + \dots + E_{-m},\quad
        E_{-k}:\Delta^n_r\to\wp_{-k}  
\end{equation}
and show  that starting from $E_{-1}  = \Gamma_{-1}$ it is possible to
construct each function $E_{-k}$ inductively using $(6.10)$. 
\vskip 10pt

\par Now, because the \lq\lq untwisted\rq\rq period map 
$\psi:\Delta^n\to\check\M$   depends upon the   choice  of coordinates
$(s_1,\dots,s_n)$ via the connection $\nabla^c$, in order to state the
next two  results we  must fix  a  choice of  holomorphic  coordinates
$(s_1,\dots,s_n)$   on $\Delta^n$   relative   to which   the  divisor
$D=\Delta^n-\Delta^{*n}$ assumes  the form  $D=\{\,p\in\Delta^n   \mid
s_1(p)\cdots s_n(p) = 0\,\}$. 

\setcounter{theorem}{18}			
\begin{theorem} An admissible VGPMHS $\V\to\Delta^{*n}$ may be recovered 
from the following data:
\begin{description}
\item{(1)} The flat connection $\nabla$ of $\V$.
\item{(2)} The Higgs field $\pd X_{-1}$.
\item{(3)} The limiting IMHM $(6.3)$.
\end{description}
\end{theorem}
\begin{proof} Since the coordinates $(s_1,\dots,s_n)$ are known, we may
recover the connection $\nabla^c$ via $(1)$.  Therefore, by using the
data contained in $(3)$, we may recover rational structure of $\V$, 
the weight filtration $\W$, and the bilinear forms ${\calS_k}$.
To recover the map $\psi:\Delta^n\to\check\M$, and hence the Hodge
filtration $\F$, observe that by holomorphicity it will suffice to 
determine $\psi$ over any open subset of $\Delta^n$.  However, by 
Theorem $(6.16)$ we can recover $\psi$ on an neighborhood of zero
using $(2)$ and $(3)$.
\end{proof}

\setcounter{corollary}{19}			
\begin{corollary} The machinery of Theorem $(6.19)$ establishes
an equivalence of categories between admissible VGPMHS $\V\to\Delta^{*n}$
which are unipotent in the sense of Definition $(5.26)$ and the 
corresponding data $(1)$--$(3)$.
\end{corollary}
\begin{proof} The key step in establishing an equivalence of categories
is that we must not shrink the domain of the Higgs field $\pd X_{-1}$.
Equivalently, the representation $\psi = e^{\Gamma}.F_{\infty}$ must
hold over all of $\Delta^n$.  To prove this, observe that by unipotency
$$
        \psi(\Delta^n) \subseteq \exp(Lie_{-1}).F_{\infty} 
$$
\end{proof}

\setcounter{example}{20}				
\begin{example} Let $\V\to\Delta^{*n}$ be an admissible variation
of Hodge--Tate type.  Then, because the variation is Hodge--Tate, the
relative weight filtration $\rel W$ coincides with the original weight
filtration $W$ and
$$
        q_{\infty} = \lam_{(F_{\infty},W)} 
$$
Moreover, upon selecting suitable branches of $\log s$, the subbundle
$\U^p$ of $\V$ may be identified with 
$$
       \U^p(s) = e^{\frac{1}{2\pi i}\sum_{j=1}^n\,\log s_j\otimes N_j}\,
                 e^{\Gamma(s)}.U^p_{\infty} 
$$

Therefore, by Theorem $(5.1)$:
\begin{eqnarray*}
  \theta &=& [\pd(e^{\frac{1}{2\pi i}\sum_{j=1}^n\, \log s_j\otimes N_j}\,
                   e^{\Gamma})]
                  [e^{\frac{1}{2\pi i}\sum_{j=1}^n\, \log s_j\otimes N_j}\,
                   e^{\Gamma}]^{-1} \\
     &=& \Omega + Ad(e^{\frac{1}{2\pi i}\sum_{j=1}^n\,\log s_j\otimes N_j})
           ((\pd e^{\Gamma}) e^{-\Gamma}) \\
     &=& \Omega + \pd\Gamma_{-1} \mod \bigoplus_{a\leq -2}\,\wp_a 
      = \pd X_{-1} \mod \bigoplus_{a\leq -2}\,\wp_a 
\end{eqnarray*}
Thus, the horizontal component of $\theta$ relative to $(F_{\infty},W)$
coincides with the Higgs field $\pd X_{-1}$.  
\end{example}

\begin{myremark} Because of the graded structure of $q_{\infty}$, there
exist universal Lie polynomials $P_2, P_3,\dots$ such that
$$
        \pd\Gamma_{-k} = P_k(\Gamma_{-1},\dots,\Gamma_{1-k},
                           \pd\Gamma_{-1},\dots,\pd\Gamma_{1-k},\Omega)
$$
whenever $\Gamma$ is a solution to $(6.10)$.  In particular, as they
shall be of use in \S 7, let us observe that the first two of the equations
are given by 
\setcounter{equation}{21}
\begin{eqnarray}
     \pd\Gamma_{-2} &=& [\Gamma_{-1},\Omega + \half\pd\Gamma_{-1}]
                        \nonumber \\
     \pd\Gamma_{-3} &=& [\Gamma_{-2},\Omega +\half\pd\Gamma_{-1}]
                     + \frac{1}{12}[\Gamma_{-1},[\Gamma_{-1},\pd\Gamma_{-1}]]
\end{eqnarray}
\end{myremark}        

\section{A--model}
\label{sec:7}
\setcounter{equation}{0}

The relationship between the Higgs fields constructed in \S 5--6 and the 
quantum cohomology of Calabi--Yau threefolds may be roughly summarized 
as follows:
\begin{description}
\item{$\bullet$} \hph{aaa} The A--model correlation functions of a Calabi--Yau 
    threefold $X$ are completely determined by the Higgs field $\pd X_{-1}$ of
    the corresponding A--model variation of Hodge structure described in 
    \cite{Cox--Katz}.

\item{$\bullet$} \hph{aaa} The B--model correlation functions of the 
    corresponding mirror family $X^{\circ}$ are completely determined by the 
    Higgs field $\theta$ of the corresponding Hodge--Tate variation
    described in \cite{Deligne}.
\end{description}
In this section, we shall treat the A--model side of the story from the 
standpoint of the Higgs field $\pd X_{-1}$.  The details of the B--model 
side will be treated in \S 8.

\par As a prelude to our discussion of the A--model variation of Hodge
structure and its relationship to quantum cohomology, we need to first
recall some standard terminology:  Let $X$ be a smooth Calabi--Yau threefold.
Then, the corresponding K\"ahler cone $K(X)$ is just the set:
\begin{equation}
	K(X) = \{\, \omega\in H^2(X,\R) \mid 
		    \omega\hph{a}\text{is K\"ahler}\,\}	
\end{equation}
Likewise, complexified K\"ahler space $K_{\C}(X)$ is just the quotient
of the set
$$
	\{\, \omega\in H^2(X,\C) \mid \Im(\omega)\in K(X)\,\} 
$$
by the torsion free part of $H^2(X,\Z)$.

\setcounter{definition}{1}				
\begin{definition} A simplicial cone $\sigma\subset H^2(X,\R)$ is said to
be a large radius limit point of $K_{\C}(X)$ provided that it is of 
maximal dimension and satisfies $Int(\sigma)\subseteq K(X)$.
\end{definition}

\par To construct a nice partial compactification of $K_{\C}(X)$ about a
given large radius limit point $\sigma$, let us suppose $\sigma$ to be
generated by a collection of basis vectors
\setcounter{equation}{2}
\begin{equation}
	T_1,\dots,T_n\in (H^2(X,\Z)/\text{torsion})\cap \overline{K(X)}	
\end{equation}
Now, by virtue of $(7.3)$
$$
	Int(\sigma) 
	 = \{\, a_1 T_1 + \cdots + a_n T_n \mid a_1,\dots, a_n > 0  \,\}
$$
and hence the product of upper half-planes
$$
	U^n_{\sigma} = \{\, u_1 T_1 + \dots + u_n T_n \mid 
				\Im(u_1),\dots,\Im(u_n) > 0 \,\}
$$
projects onto a neighborhood of infinity $D_{\sigma}$ of $K_{\C}(X)$ via
the quotient map
\begin{equation}
  	\vec u\in U^n_{\sigma}\mapsto [\vec u] \in K_{\C}(X)
\end{equation}
Consequently, we may use the quotient map $(7.4)$ to obtain a completion
$$
	D_{\sigma}\stackrel{i}{\hookrightarrow}\overline{D}_{\sigma}
$$
via the standard covering map
$$
	(u_1,\dots,u_n)\in U^n_{\sigma} \mapsto 
	(e^{2\pi i u_1},\dots,e^{2\pi i u_n})\in\Delta^{*n}
$$
and the inclusion $\Delta^{*n}\stackrel{i}{\hookrightarrow}\Delta^n$.

\par In particular, upon setting $q_j = e^{2\pi i u_j}$, we obtain a natural
system of coordinates
$$
	(q_1,\dots,q_n):D_{\sigma}\to\Delta^{*n}
$$ 
relative to which the completion 
$D_{\sigma}\stackrel{i}{\hookrightarrow}\overline{D}_{\sigma}$ becomes the
standard embedding $\Delta^{*n}\stackrel{i}{\hookrightarrow}\Delta^n$.
\vskip 10pt

\par In light of these observations, we shall henceforth assume 
that our large radius limit point $\sigma$ comes equipped with a choice 
of basis $T_1,\dots, T_n$ of $H^2(X,\Z)/\text{torsion}$ relative to 
which $(7.3)$ holds. 
\vskip 10pt

\par The next ingredient the we must assemble in order to discuss the 
A--model connection is the small quantum product $a\ast b$.

\setcounter{lemma}{4}						
\begin{lemma} Let $\sigma$ be a large radius limit point of $K_{\C}(X)$.
Then, relative to the coordinates 
$$
	\omega(u_1,\dots,u_n) = \sum_{j=1}^n u_j T_j
$$
defined on $U^n_{\sigma}$ by the basis $(7.3)$, the Gromov--Witten
potential $\Phi$ of $X$ assumes the form 
$$
      \Phi(u_1,\dots,u_n) = \left(\frac{1}{6}\int_X\, \omega^3\right)
                            + \Phi_{hol}(q_1,\dots,q_n)        
$$
with respect to the formal $q$-power series:
\setcounter{equation}{5}
\begin{equation}
        \Phi_{hol}(q_1,\dots,q_n) = 
        \frac{1}{(2\pi i)^3}\sum_{\b\in H^2(X,\Z)-\{0\}}
	   \langle I_{0,0,\b}\rangle\,
	    e^{2\pi i\int_{\b}\,\omega}
\end{equation}
\end{lemma}

\setcounter{corollary}{6}				
\begin{corollary} Let $\sigma$ be a large radius limit point of $K_{\C}(X)$
and 
$$
	{}^{\vee}:H^k(X,\C)\to H^{6-k}(X,\C)
$$
denote the map induced by Poincar\'e duality.  Then,
$$
	T_j\ast T_k = \sum_{\ell}\,
		\frac{\pd^3 \Phi}{\pd u_j\pd u_k\pd u_{\ell}} T_{\ell}^{\vee}
$$
\end{corollary}

\begin{myremark} The formal $q$-power series $(7.6)$ is expected to converge
for all values of $q$ sufficiently small.  Therefore, in order to simplify
our exposition, for the remainder of this section we shall assume the sum 
$(7.6)$ to converge over all of $\overline{D}_{\sigma}$.  
\end{myremark}

\par To construct the A--model variation of Hodge structure, let 
$$	
	H = \bigoplus_p\,H^{p,p}(X)
$$ 
and $\V$ denote the product bundle $H\times D_{\sigma}$ endowed
with the corresponding Dubrovin connection:
\setcounter{equation}{7}
\begin{equation}
   \nabla = d + A,\qquad A_{\frac{\pd}{\pd u_j}}\a = T_j\ast\a    
\end{equation}
Then, as may be easily checked by direct computation, the associativity of 
the small quantum product $\ast$ implies that $\nabla$ is flat.

\par Thus, in order to finish the construction of the A--model variation, 
it remains for us to describe: 
\begin{description}
\item{(1)} The Hodge filtration $\F$.
\item{(2)} The polarization $Q$.
\item{(3)} The integral structure $\V_{\Z}$.
\end{description}

\par To endow $\V$ with an integral structure which is flat with respect to
$\nabla$, let $\V^c= H\times\overline{D}_{\sigma}$ denote the canonical 
extension of $\V$ and 
$$
	\nabla^c = \nabla 
	- \frac{1}{2\pi i}\sum_{j=1}^n\,\frac{dq_j}{q_j}\otimes N_j,
	\qquad N_j = -\log(T_j)
$$
denote the corresponding connection of $\V^c$ defined by the choice of
coordinates $q_1 = e^{2\pi i u_1},\dots, q_n = e^{2\pi i u_n}$.  Then, in 
keeping with our discussions in \S 6, it is natural to define it is natural 
to define
$$
	\sigma = e^{-\sum\,\frac{1}{2\pi i}\log(q_j)N_j}\sigma^c
$$
to be a flat (multivalued) section of $\V_{\Z}$ if and only if $\sigma^c$
is a $\nabla^c$-flat section of $\V^c$ such that 
$\sigma^c(0)\in H^{\ast}(X,\Z)$.

\par To polarize $\V$, we pair $\a\in H^{p,p}(X)$ and 
$\b\in H^{3-p,3-p}(X)$ by the rule:
$$
        Q(\a,\b) = (-1)^p\int_X\,\a\wedge\b 
$$
Likewise, following \cite{Cox--Katz}, we define:
$$
        \F^p = \bigoplus_{a\leq 3-p}\, H^{a,a}(X,\C) 
$$

\setcounter{theorem}{8}					
\begin{theorem}{\cite{Cox--Katz}} The data $(\V_{\Z},\F,Q)$ defines a 
variation of pure, polarized Hodge structure of weight $3$ over a neighborhood
of zero in $D_{\sigma}$ for which the action of the monodromy logarithm 
$N_j$ on the central fiber of $\V^c$ coincides with cup product with $T_j$.
Thus, by virtue of the Hard Lefschetz Theorem, the weight filtration of
the monodromy cone
$$
	{\cal C} = \{\,a_1 N_1 + \dots + a_n N_n \mid a_1,\dots,a_n>0\,\}
$$
is given by the rule:
$$
        W_k = \bigoplus_{2a\geq 6-k}\,H^{a,a}(X) 	
$$

\par Moreover, upon letting ${}^{\vee}:H^k(X,\C)\to H^{6-k}(X,\C)$ 
denote the map induced by Poincar\'e duality and $T_0=1\in H^0(X,\C)$, 
the following sections may be shown to define a $\nabla^c$-flat framing of 
$\V^c$:
\setcounter{equation}{9}
\begin{alignat}{2}
  &\sigma_j = T_j 
   - \sum_{\ell}\,\frac{\pd^2\Phi_{hol}}{\pd u_j\pd u_{\ell}}\,T_{\ell}^{\vee} 
   + \frac{\pd\Phi_{hol}}{\pd u_j}\,T_0^{\vee},
  & \qquad &\sigma^j = T_j^{\vee} \nonumber \\
  &\sigma_0 = T_0 
    - \sum_{\ell}\,\frac{\pd\Phi_{hol}}{\pd u_{\ell}}\,T_{\ell}^{\vee}
    + 2\Phi_{hol}\,T_0^{\vee},
  & \qquad &\sigma^0 = T_0^{\vee} 
\end{alignat}  					
\end{theorem}

\par Armed with these preliminaries, we can now rewrite the Hodge filtration 
$\F$ of $\V$ in terms of the $\nabla^c$-flat frame $(7.10)$ and thus obtain 
the \lq\lq untwisted\rq\rq{} period map $\psi:\Delta^n\to\check{\cal{D}}$
described in \S 6.  Upon carrying out these computations, one finds that
\begin{eqnarray} 
  \F^3 &=& \text{span}_{\C}(\sigma_0 
  + \sum_{\ell}\,\frac{\pd\Phi_{hol}}{\pd u_{\ell}}\,\sigma^{\ell} 
  - 2\Phi_{hol}\,\sigma^0) \nonumber \\
  \F^2 &=& \F^3\oplus\text{span}_{\C}(\sigma_j 
  + \sum_{\ell}\,\frac{\pd^2\Phi_{hol}}{\pd u_j\pd u_{\ell}}\,\sigma^{\ell}
  - \frac{\pd\Phi_{hol}}{\pd u_j}\sigma^0) \\
  \F^1 &=& \F^2\oplus\text{span}_{\C}(\sigma^j),\qquad  
  \F^0 = \F^1\oplus\text{span}_{\C}(\sigma^0)  \nonumber
\end{eqnarray}
where the index $j$ ranges from 1 to $n$.  

\par Now, as may be easily checked, relative to definitions of $\F$ and $\W$
given above:
\begin{equation}
        I^{p,p}_{(F_{\infty},W)} = H^{3-p,3-p}(X,\C)	
\end{equation}
Therefore, by virtue of equation $(7.11)$, 
$$
	\psi(q_1,\dots,q_n) = e^{\Gamma(q_1,\dots,q_n)}.F_{\infty}
$$
with $\Gamma = \Gamma_{-1} + \Gamma_{-2} + \Gamma_{-3}$, and
\begin{alignat}{2}
          &\Gamma_{-1}(T_k) = \sum_{\ell=1}^n\,
            \frac{\pd^2 \Phi_{hol}}{\pd u_k\pd u_{\ell}}\,T_{\ell}^{\vee} 
          & \qquad 
          & \Gamma_{-2}(T_k) = -\frac{\pd \Phi_{hol}}{\pd u_k}\,T_0^{\vee} 
	    \nonumber \\
          &\Gamma_{-2}(T_0) 
            = \sum_{\ell}\,\frac{\pd\Phi_{hol}}{\pd u_{\ell}}\,T_{\ell}^{\vee} 
          & \qquad 
	  &\Gamma_{-3}(T_0) = -2\Phi_{hol}\,T_0^{\vee}
\end{alignat}

\setcounter{theorem}{13}				
\begin{theorem} Let $\sigma$ be a large radius limit point of $K_{\C}(X)$ and 
$$
	\pd X_{-1} 
	= \sum_{j=1}^n \left(N_j + \frac{\pd\Gamma_{-1}}{\pd u_j}\right)
		       \otimes du_j
$$
be the Higgs field of the corresponding A--model variation $\V$.
Then, for any element $\a\in H$:
$$
        \pd X_{-1}\left(\frac{\pd}{\pd u_j}\right)\,\a = T_j\ast \a
$$
\end{theorem}
\begin{proof}  Since $N_j$ acts as cup product by $T_j$ and 
$T_j\ast\a = T_j\cup\a$ provided $\a$ has no component in $H^2(X,\C)$,
it will suffice to establish the claim for $\a\in H^2(X,\C)$.  However,
by equation $(7.13)$, 
$$
        \pd X_{-1}\left(\frac{\pd}{\pd u_j}\right)\,T_k =
        T_j\wedge T_k + 
        \sum_{\ell}\,\frac{\pd^3\Phi_{hol}}{\pd u_j\pd u_k\pd u_{\ell}}\,
        T_{\ell}^{\vee}  
        = T_j\ast T_k 			
$$
\end{proof}

\setcounter{definition}{14}				
\begin{definition} Let $M$ be a smooth projective variety, $e_0,\dots,e_m$
be a basis for the rational cohomology of $M$ each term of which is 
homogeneous with respect to the degree map
$$
	deg:H^*(M,\Q)\to\Z
$$
and $t_0,\dots,t_m$ be the corresponding supercommuting variables defined
by the rule:
$$
	t_j t_k = (-1)^{deg(e_j) deg(e_k)} t_k t_j
$$
Then, a formal power series $\Phi(t_0,\dots,t_m)$ is said to be solution of 
the WDVV (Witten--Dijgraaf--Verlinde--Verlinde) equation provided the 
formal product
\setcounter{equation}{15}
\begin{equation}
	e_j \ast e_k 
	= \sum_{\ell=0}^m\, 
	  \left(\frac{\pd^3\Phi}{\pd t_j\pd t_k\pd t_{\ell}}\right) 
	   e_{\ell}^{\vee}
\end{equation}
turns $H^*(M,\C)$ into a supercommutative ring with identity element
$1 \in H^0(X,\C)$.
\end{definition}

\setcounter{example}{16}				
\begin{example} The standard cup product of $M$ defines a natural solution
to the WDVV equations via the potential function
$$
	\Phi_0(t_0,\dots,t_m) 
	  = \int_M\,\left(\sum_{j=0}^m\,t_j e_j\right)^3
$$
\end{example}

\par Another example of a solution to the WDVV equation is provided by the
small quantum product.  To relate these two solutions, observe that in the
case of a Calabi--Yau threefold $X$, the small quantum product coincides
with the usual cup product except on elements of $H^2(X,\C)$.  Moreover, 
upon selecting a large radius limit point $\sigma$ of $K_{\C}(X)$ and setting 
$e_j = T_j$ for  $j=1,\dots,n$, it follows by Corollary $(7.7)$ that
\setcounter{equation}{17}
\begin{equation}
	\Phi(t_0,\dots,t_m) = \Phi_0(t_0,\dots,t_m) 
			   + \Phi_{hol}(q_1,\dots,q_n)		
\end{equation}
is a potential function for the small quantum product.  In addition, because
\begin{equation}
	\Phi_{hol}(0,\dots,0) = 0				
\end{equation}
the small quantum product is in fact a deformation of the standard potential
function $\Phi_0$.

\setcounter{theorem}{19}				     
\begin{theorem} Let $\sigma$ be a large radius limit point of the Calabi--Yau
threefold $X$.  Then, there exists a bijective correspondence between germs 
of holomorphic solutions of the WDVV equations which satisfy conditions 
$(7.18)$--$(7.19)$ and the germs of variations of pure, polarized Hodge 
structures of weight $3$ on $D_{\sigma}$ which are governed by equation 
$(7.13)$.
\end{theorem}
\begin{proof} Given a variation of Hodge structure $\V$ governed by equation 
$(7.13)$, one defines the corresponding potential $\Phi$ by the rule:
\setcounter{equation}{20}
\begin{equation}
	\Phi(t_0,\dots,t_m) = \Phi_0(t_0,\dots,t_m) 
				+ \Phi_{hol}(q_1,\dots,q_n)	
\end{equation}
where
$$
	\Gamma_{-3}(T_0) = -2\Phi_{hol}\,T_0^{\vee}
$$

\par Now, in order to verify that equation $(7.21)$ defines a solution  of
the WDVV equation, observe that because $X$ is a Calabi--Yau threefold, it
will suffice to check that
$$
	T_a\ast(T_b\ast T_c) = (T_a\ast T_b)\ast T_c
$$
for all indices $a$, $b$ and $c$.  

\par To this end, note that by virtue of our assumption that $e_j = T_j$
for $j=1,\dots,n$ and equation $(7.13)$ we have:
$$
	T_a\ast(T_b\ast T_c) 
	= \pd X_{-1}(\frac{\pd}{\pd t_a})\circ\pd X_{-1}(\frac{\pd}{\pd t_b})
	  \circ\pd X_{-1}(\frac{\pd}{\pd t_c})\,T_0
$$
Therefore, on account of the symmetry condition 
$\pd X_{-1}\wedge\pd X_{-1} = 0$,
\begin{eqnarray*}
	T_a\ast(T_b\ast T_c) 
	&=& \pd X_{-1}(\frac{\pd}{\pd t_a})\circ\pd X_{-1}(\frac{\pd}{\pd t_b})
	  \circ\pd X_{-1}(\frac{\pd}{\pd t_c})\,T_0 \\
	&=& \pd X_{-1}(\frac{\pd}{\pd t_c})\circ\pd X_{-1}(\frac{\pd}{\pd t_a})
	  \circ\pd X_{-1}(\frac{\pd}{\pd t_b})\,T_0 \\
	&=& T_c\cup (T_a\ast T_b) = (T_a\ast T_b)\cup T_c 
	    = (T_a\ast T_b)\ast T_c
\end{eqnarray*}
since the binary operation $\ast$ defined by equation $(7.22)$ agrees with
cup product except on classes of degree 2.

\par To verify that this construction is compatible with passage to the germ
of $\V$, observe that if $\tilde\V$ is variation of Hodge structure which is
governed by equation $(7.13)$  and agrees with $\V$ over some neighborhood
$O$ of zero in $D_{\sigma}$ then 
$$
	\Phi|_O = \tilde\Phi|_O
$$

\par To establish the converse, let $\Phi$ be a holomorphic solution of the
WDVV equation which satisfies conditions $(7.18)$ and $(7.19)$.  Then, 
because the potential function $\Phi$ is an honest holomorphic function
and not merely a formal power series, we may use conditions $(7.18)$ and
$(7.19)$ to define a candidate Hodge filtration $\F$ over some neighborhood
of zero  in $D_{\sigma}$ via equation $(7.13)$.  

\par Likewise, we may use the Dubrovin connection $\nabla$ defined by equation
$(7.8)$ and our potential function $\Phi$ to define a candidate for the 
integral structure $\V_{\Z}$ and the polarization $Q$.

\par Now, as the reader will see upon reviewing the results of \S 6, in order
to verify that the triple $(\V_{\Z},\F,Q)$ defines a variation of pure,
polarized Hodge structure of weight $3$ over some neighborhood of zero in
$D_{\sigma}$, it will suffice to check the following three conditions:
\begin{description}
\item{(1)} \hph{aa} The Higgs field $\pd X_{-1} = \Omega + \pd\Gamma_{-1}$ 
    defined by equation $(7.13)$ satisfies the symmetry condition 
    $\pd X_{-1}\wedge\pd X_{-1} = 0$.

\item{(2)} \hph{aa} The function $\Gamma =\Gamma_{-1}+\Gamma_{-2}+\Gamma_{-3}$ 
    constructed by Theorem $(6.16)$ from the Higgs field $\pd X_{-1}$ is 
    compatible with the given definitions of $\Gamma_{-2}$ and $\Gamma_{-3}$ 
    appearing in equation $(7.13)$.
\item{(3)} \hph{a} The triple $(\V_{\Z},\F,Q)$ defines a polarized Hodge 
    structure at each point $q$ of some neighborhood of zero in $D_{\sigma}$.
\end{description}

\par To verify condition $(1)$, one first recalls the commutativity and 
associativity properties of the product $\ast$ and then simply computes
the action of the two form $\pd X_{-1}\wedge\pd X_{-1}$ on $\V$.

\par To verify condition $(2)$, one simply checks that the functions 
$\Gamma_{-2}$ and $\Gamma_{-3}$ defined by equation $(7.13)$ are a 
solution of the system of differential equations
\begin{eqnarray*}
     \pd\Gamma_{-2} &=& [\Gamma_{-1},\Omega + \half\pd\Gamma_{-1}] 	\\
     \pd\Gamma_{-3} &=& [\Gamma_{-2},\Omega +\half\pd\Gamma_{-1}]
                     + \frac{1}{12}[\Gamma_{-1},[\Gamma_{-1},\pd\Gamma_{-1}]]
\end{eqnarray*}
derived at the end of \S 6.

\par Finally, to verify condition $(3)$, one may proceed as in \cite{Cox--Katz}
and simply check that the limiting mixed Hodge structure $(F_{\infty},W)$ of
$\V$ is polarized by the bilinear form $Q$ and the monodromy cone
$$
	{\cal C} = \{\, a_1 N_1 + \dots + a_n N_n \mid a_1,\dots,a_n > 0\,\}
$$

\par Now, in order to prove that the germ of $\V$ only depends upon the germ
of $\Phi$, observe that if $\tilde\Phi$ is a potential function which agrees
with $\Phi$ over some neighborhood $O$ of zero in $D_{\sigma}$ then
$$
	\V|_O = \tilde\V|_O
$$
\end{proof}

\begin{myremark} In order to avoid having to work at the level of germs in the
preceding theorem, observe that equation $(7.13)$ also governs a Hodge--Tate
variation since the Hodge filtration $\F$ is opposed to the monodromy weight 
filtration $\W$ over all of $D_{\sigma}$. 
\end{myremark}

\par Finally a reminder:  All of the results stated in this section are 
predicated upon the implicit assumption that our given large radius 
limit point $\sigma$ satisfies the integrality condition $(7.3)$.
\section{B model}
\label{sec:8}
\setcounter{equation}{0}

We shall now revisit \cite{Deligne} from the viewpoint of the Higgs field 
carried by a Hodge--Tate variation.  For reference, the rough outline of 
\cite{Deligne} is as follows:
\begin{description}
\item{(1)} \hph{aa} In the vicinity of a maximally unipotent boundary point, 
	   the variations of pure, polarized Hodge structure $\cal H$ arising 
	   in the B--model of mirror symmetry give rise to a corresponding 
	   variation of mixed Hodge structure $\V\to\Delta^{*n}$ which is
  	   obtained by pairing the Hodge filtration $\F$ of $\cal H$ with the
	   monodromy weight filtration $\W$.

\item{(2)} \hph{aa} The variations of mixed Hodge structure described in (1) 
	   are of Hodge-Tate type, defined over $\Z$, and may therefore 
	   described in terms of extension classes (sections)
           \begin{equation}
              E_{\ast}\in 
              End_{\Z}(Gr^{\W})_{-2}\otimes{\cal{O}}^*_{mer}(\Delta^{*n})
           \end{equation}
           Moreover, given a generator $1\in Gr^W_0(\Z)$,
           \begin{equation}
                Gr^W_{-2} = \text{span}_{\Z}(N_1(1),\dots,N_n(1)) 
           \end{equation}
\item{(3)} \hph{aa} By virtue of $(8.1)$ and $(8.2)$, there exists canonical
           coordinates $(q_1,\dots,q_n)$ on $\Delta^n$ relative to
           which $E_0:Gr^W(\Z)\to Gr^W_{-2}(\Z)$ assumes the form 
           \begin{equation}
               E_0(1) = \sum_{j=1}^n\, q_j N_j(1)  
           \end{equation}
           
\item{(4)} \hph{aa} The Yukawa coupling of the underlying physical theory
           may be computed as the \lq\lq logarithmic derivative\rq\rq{} 
           of  $E_1:Gr^{\W}_{-2}\to Gr^{\W}_{-4}$.
\end{description}

\par Now, for Hodge-Tate variation $\V$, the corresponding Deligne grading 
$\Y$ determines both the Hodge filtration and weight filtration of $\V$.  
Moreover the connection $\nabla$ of $\V$ may be expressed as
$$
        \nabla = \nabla^h + \theta       
$$
with $\nabla^h = \bar\pd + \pd$ preserving each subbundle $I^{p,p}$ of 
$\V$ and
\begin{equation}
        \theta\in \Omega^1\otimes \gg^{-1,-1} 
\end{equation}

\par Therefore, as discussed in \S 8 of \cite{Deligne}, the vector bundle 
isomorphism $\Y:\V\cong Gr^{\W}$ allows the transport of the flat connection 
$\nabla$ of $\V$ to the connection
\begin{equation}
        D = d -\frac{1}{2\pi i}\log E    
\end{equation}
on $Gr^{\W}$, where $E$ is the extension class of $\V$ and 
$d$ denotes the standard connection on $Gr^{\W}$ defined by the 
integral structure $Gr^{\W}(\Z)$.

\par Unraveling definitions, direct computation shows that the connection
appearing in $(8.5)$ may be rewritten as
$$
        D = d + \theta   
$$
relative to the natural action 
\begin{equation}
    \theta:T^{hol}(\Delta^{*n})\otimes Gr^{\W}_{2p}\to 
                Gr^{\W}_{2p-2} 
\end{equation}
induced on each $Gr^{\W}_{2p}$ by $(8.4)$.  Comparing $(8.5)$ and 
$(8.6)$ it therefore follows that:
\begin{equation}
        \theta = -\frac{1}{2\pi i}d\log E 
\end{equation}
In particular, the Yukawa coupling is given by the action 
$$
    d\log E_{-2} = -2\pi i\theta\,:\,
    T^{hol}(\Delta^{*n})\otimes Gr^{\W}_{-2}\to Gr^{\W}_{-4} 
$$
\begin{myremark} In \cite{Deligne}, $\V$ is a variation in homology, hence 
$Gr^{\W}_k = 0$ for $k>0$.  While I shall keep Deligne's notation, I will 
continue to refer to $\V$ in terms of cohomology.
\end{myremark}

\par Suppose now that the fibers of the constant local system $Gr^{\W}$ 
arise as the middle cohomology of an algebraic variety (the mirror), or
more generally a family of such varieties.  Then, there must be a natural
cup product structure on $Gr^{\W}$ as well as a notion of Poincar\'e
duality (assuming the mirror variety has only mild singularities).  As 
discussed by \cite{Deligne}, these structures should be induced by the
pairing\footnote{Note: $Q$ is the pairing on $Gr^{\W}$ 
obtained from the $SL_2$ orbit theorem.  See \cite{CK} for details.} 
$Q:Gr^W_{-k}(\Z)\otimes Gr^W_{-6-k}(\Z)\to\Z$ determined by the variation 
$\V$ and the algebra structure $Gr^{\W}(\Z)$ inherits from the 
polynomial algebra $\Z[N_1,\dots, N_n]$. 

\par More precisely, let
$$
  L = \text{span}_{\Z}(N_1,\dots,N_n)   
$$
Then, by virtue of $(8.2)$, we obtain an isomorphism 
$$
  Gr^W \cong Z\oplus L \oplus L^{\vee}\oplus Z^{\vee}   
$$
via the maps
$$
\begin{array}{l}
       k\in\Z\mapsto k(1)\in Gr^W_0(\Z)  \\
       N\in L\mapsto N(1)\in Gr^W_{-2}(\Z) 
\end{array}				
$$
and the pairing $Q$.  Under this isomorphism, the \lq\lq cup product\rq\rq{}
map $L\otimes L\to L^{\vee}$ is determined by the functional
$$
        (xy)(z) := Q(1,xyz(1))     
$$

\par To obtain a deformation of this product structure on $Gr^{\W}$,
observe that by virtue of $(6.21)$, 
$$
        \Res_{q_j = 0}(\theta) = \frac{1}{2\pi i} N_j    
$$
and hence introducing the vector fields $\xi_j = 2\pi iq_j\frac{\pd}{\pd q_j}$
it follows that 
$$
        \theta(\xi_j) 
          = N_j + (\text{higher order terms}) 
$$
Consequently, the product structure defined by the three point function
\begin{equation}
        \phi(N_a, N_b, N_c) 
        = Q(1,\theta(\xi_a)\circ\theta(\xi_b)\circ\theta(\xi_c)\,1)
\end{equation}
will be a deformation of the monodromy algebra of $Gr^{\W}(\Z)$.

\setcounter{theorem}{8}			         
\begin{theorem} The product structure determined by $(8.8)$ 
is commutative and associative.
\end{theorem}
\begin{proof} Let $\a\cup\b$ denote the standard cup product on 
$Gr^{\W}$ defined previously and set
\setcounter{equation}{9}
$$
     \a\ast\b = 
           \left\{\begin{array}{ll}
             \a\cdot\b   &\qquad\text{if}\hph{a} \a,\b\in Gr^{\W}_{-2} \\
             \a\cup\b    &\qquad\text{otherwise}
	     \end{array}\right.
$$
where $\a\cdot\b$ is defined by the requirement
\begin{equation}
   (N_a\cdot N_b)\cup N_c = \phi(\xi_a,\xi_b,\xi_c)
\end{equation}
Therefore, by virtue of the symmetry condition $\theta\wedge\theta = 0$, the
value of $(8.10)$ is invariant under the natural action of the permutation
group $S_3$ on the labels $a, b, c$.  Consequently,
$$
        (N_a\cdot N_b)\cup N_c = (N_b\cdot N_c)\cup N_a = N_a\cup(N_b\cdot N_c)
$$
on account of $S_3$--invariance and the commutativity of the product between 
$Gr^{\W}_{-2}$ and $Gr^{\W}_{-4}$.  Likewise, 
$N_a\cdot N_b = N_b\cdot N_a$ because by $S_3$--invariance
$$
        (N_a\cdot N_b)\cup N_c = (N_b\cdot N_a)\cup N_c
$$
\end{proof}

\begin{myremark} It is not necessary to explicitly determine the coordinates
$(q_1,\dots,q_n)$ in order to form the quantum product via this procedure.
More precisely, we only need the vector fields 
$\xi_1,\dots, \xi_n$ which may be recovered as follows: The map
$$
        \theta:T^{hol}(\Delta^{*n})\otimes Gr^{\W}_0\to Gr^{\W}_{-2}
$$
must have an expression of the form
$$
        \theta(\xi)\,1 = \sum_j\, N_j(1)\otimes\Omega_j(\xi)   
$$
relative to meromorphic 1--forms $\Omega_1,\dots,\Omega_n$.  However, 
comparing $(8.3)$ and $(8.7)$, it follows that $\{\xi_j\}$ is
the dual frame of $\{\Omega_j\}$ up to a factor of $(-1)$.
\end{myremark}

\appendix
\section{Appendix}


\par In section I communicate an alternative proof of Theorem $(5.1)$
due to P. Deligne.  With the exception of a few minor adjustments and remarks,
our presentation essentially follows \cite{Deligne2} verbatim.
\vskip 10pt

\par For any mixed Hodge structure $(F,W)$, one has\newline
\hph{aaaaaa}\newline
$(a)$ $F^p = \oplus_{r\geq p}\, I^{r,s}$.\newline
\hph{aaaaaa}\newline
$(b)$ Let $\bar I$ be obtained by interchanging the roles of $F$ 
and $\bar F$ in the definition of $I$, i.e.
$$
	(\bar I)^{r,s} := \overline{I^{s,r}}
$$
Then,
\begin{equation}
	\bigoplus_{s\leq q} (\bar I)^{r,s} = \bigoplus_k\, W_k\cap F^{k-q}
\end{equation}
\vskip 10pt

\par Indeed, as the preceding sum, together with a bi-index $(r,s)$
contains all $(r^{\prime},s^{\prime}) < (r,s)$, one may replace $\bar I$
by  $I$ in $(A.1)$.  To prove that
\begin{equation}
	\bigoplus_{s\leq q}\, I^{r,s} = \bigoplus_k\, W_k\cap F^{k-q}
\end{equation}
observe that
\begin{description}
\item{(1)} \hph{aa} By definition,
$$
	   W_k\cap F^{k-q} = \bigoplus_{r+s\leq k,\,r\geq k-q}\, I^{r,s}
$$
In particular, $W_k\cap F^{k-q}\subseteq \bigoplus_{s\leq q}\, I^{r,s}$
since 
$$
	r+s\leq k,\hphantom{a} r\geq k-q \implies r+s\leq k \leq r+q 
	\implies s\leq q
$$
\item{(2)} \hph{aa} Conversely, 
$$
	 \bigoplus_{s\leq q}\, I^{r,s} \subseteq\bigoplus_k\, W_k\cap F^{k-q}
$$
since given $(r,s)$ with $s\leq q$ there exists an integer $k$ such that
$$
	r+s \leq k \leq r+q
$$
(i.e. $r+s\leq k$ and $r\geq k-q$).
\end{description}  

\par Now, if $\Phi_{\ast}$ is the increasing filtration $(b)$ attached to a
variation of graded-polarized mixed Hodge structure $\V$ then Griffiths
transversality for for the Hodge filtration $\F$ gives rise to a 
\lq\lq transversality\rq\rq{} condition
\begin{equation}
	\nabla\Phi_q \subseteq \Omega^1\otimes\Phi_{q+1}     
\end{equation}

\par For the $C^{\infty}$ bundle $E$ underlying such a variation, $\F$ and 
$\bar\Phi$ are the two filtrations attached to the decomposition
\begin{equation}
	E = \bigoplus_p\, \U^p,\qquad \U^p = \bigoplus_q\,\I^{p,q}
\end{equation}
(i.e. $\F^p = \bigoplus_{a\geq p}\, \U^a$ and 
$\overline{\Phi}_q = \bigoplus_{a\leq q}\, \U^a$.)  That $\F$ and $\Phi$ are
holomorphic and obey transversality is equivalent to the fact that the
direct sum decomposition $(A.4)$ defines a complex variation of Hodge
structure with respect to $\nabla$.

\begin{myremark} In addition to its directness, this approach has the
following advantages over the proof of Theorem $(5.1)$ presented
in \S 5:
\begin{description}
\item{$\bullet$} \hph{aa} It avoids the introduction of the classifying spaces 
of graded-polarized mixed Hodge structures $\M$ constructed in \S 4.
\item{$\bullet$} \hph{aa} It shows that we may drop the the assumption of 
graded-polarizability from the statement of Theorem $(5.1)$.
\end{description}
\end{myremark}

\end{document}